\documentclass{amsart}

\usepackage{amssymb}
\usepackage{amsmath}
\usepackage{array}
\usepackage{graphicx}
\usepackage{multirow}
\usepackage{xcolor}

\usepackage{mathtools}

\newcommand{\Pf}{{\em Proof}. }
\newcommand{\EPf}{\hfill$\square$}

\newcommand{\Lg}{\mbox{$\mathfrak g$}}

\newcommand{\Lk}{\mbox{$\mathfrak k$}}

\newcommand{\Lt}{\mbox{$\mathfrak t$}}

\newcommand{\x}{\mbox{$\frac12$}}
\newcommand{\codim}{\mbox{$\mathop{\mathrm{codim}}$}}

\newcommand\liegr{\sf}

\newcommand{\SU}[1]{\mbox{${\liegr SU}(#1)$}}

\newcommand{\U}[1]{\mbox{${\liegr U}(#1)$}}
\newcommand{\SP}[1]{\mbox{${\liegr Sp}(#1)$}}
\newcommand{\SO}[1]{\mbox{${\liegr SO}(#1)$}}

\newcommand{\OG}[1]{\mbox{${\liegr O}(#1)$}}
\newcommand{\Spin}[1]{\mbox{${\liegr Spin}(#1)$}}
\newcommand{\G}{\mbox{${\liegr G}_2$}}
\newcommand{\F}{\mbox{${\liegr F}_4$}}
\newcommand{\E}[1]{\mbox{${\liegr E}_{#1}$}}

\newcommand\fieldsetc{\mathbb}

\newcommand{\Z}{\fieldsetc{Z}}
\newcommand{\R}{\fieldsetc{R}}
\newcommand{\C}{\fieldsetc{C}}
\newcommand{\Q}{\fieldsetc{H}}

\newtheorem{thm}{Theorem}[section]

\newtheorem{lem}[thm]{Lemma}

\theoremstyle{remark}
\newtheorem{rmk}[thm]{Remark}

\begingroup\makeatletter\ifx\SetFigFont\undefined
\def\x#1#2#3#4#5#6#7\relax{\def\x{#1#2#3#4#5#6}}%
\expandafter\x\fmtname xxxxxx\relax \def\y{splain}%
\ifx\x\y   
\gdef\SetFigFont#1#2#3{%
  \ifnum #1<17\tiny\else \ifnum #1<20\small\else
  \ifnum #1<24\normalsize\else \ifnum #1<29\large\else
  \ifnum #1<34\Large\else \ifnum #1<41\LARGE\else
     \huge\fi\fi\fi\fi\fi\fi
  \csname #3\endcsname}%
\else
\gdef\SetFigFont#1#2#3{\begingroup
  \count@#1\relax \ifnum 25<\count@\count@25\fi
  \def\x{\endgroup\@setsize\SetFigFont{#2pt}}%
  \expandafter\x
    \csname \romannumeral\the\count@ pt\expandafter\endcsname
    \csname @\romannumeral\the\count@ pt\endcsname
  \csname #3\endcsname}%
\fi
\fi\endgroup

\begingroup\makeatletter\ifx\SetFigFont\undefined
\def\x#1#2#3#4#5#6#7\relax{\def\x{#1#2#3#4#5#6}}%
\expandafter\x\fmtname xxxxxx\relax \def\y{splain}%
\ifx\x\y   
\gdef\SetFigFont#1#2#3{%
  \ifnum #1<17\tiny\else \ifnum #1<20\small\else
  \ifnum #1<24\normalsize\else \ifnum #1<29\large\else
  \ifnum #1<34\Large\else \ifnum #1<41\LARGE\else
     \huge\fi\fi\fi\fi\fi\fi
  \csname #3\endcsname}%
\else
\gdef\SetFigFont#1#2#3{\begingroup
  \count@#1\relax \ifnum 25<\count@\count@25\fi
  \def\x{\endgroup\@setsize\SetFigFont{#2pt}}%
  \expandafter\x
    \csname \romannumeral\the\count@ pt\expandafter\endcsname
    \csname @\romannumeral\the\count@ pt\endcsname
  \csname #3\endcsname}%
\fi
\fi\endgroup

\begingroup\makeatletter\ifx\SetFigFont\undefined
\def\x#1#2#3#4#5#6#7\relax{\def\x{#1#2#3#4#5#6}}%
\expandafter\x\fmtname xxxxxx\relax \def\y{splain}%
\ifx\x\y   
\gdef\SetFigFont#1#2#3{%
  \ifnum #1<17\tiny\else \ifnum #1<20\small\else
  \ifnum #1<24\normalsize\else \ifnum #1<29\large\else
  \ifnum #1<34\Large\else \ifnum #1<41\LARGE\else
     \huge\fi\fi\fi\fi\fi\fi
  \csname #3\endcsname}%
\else
\gdef\SetFigFont#1#2#3{\begingroup
  \count@#1\relax \ifnum 25<\count@\count@25\fi
  \def\x{\endgroup\@setsize\SetFigFont{#2pt}}%
  \expandafter\x
    \csname \romannumeral\the\count@ pt\expandafter\endcsname
    \csname @\romannumeral\the\count@ pt\endcsname
  \csname #3\endcsname}%
\fi
\fi\endgroup

\begingroup\makeatletter\ifx\SetFigFont\undefined
\def\x#1#2#3#4#5#6#7\relax{\def\x{#1#2#3#4#5#6}}%
\expandafter\x\fmtname xxxxxx\relax \def\y{splain}%
\ifx\x\y   
\gdef\SetFigFont#1#2#3{%
  \ifnum #1<17\tiny\else \ifnum #1<20\small\else
  \ifnum #1<24\normalsize\else \ifnum #1<29\large\else
  \ifnum #1<34\Large\else \ifnum #1<41\LARGE\else
     \huge\fi\fi\fi\fi\fi\fi
  \csname #3\endcsname}%
\else
\gdef\SetFigFont#1#2#3{\begingroup
  \count@#1\relax \ifnum 25<\count@\count@25\fi
  \def\x{\endgroup\@setsize\SetFigFont{#2pt}}%
  \expandafter\x
    \csname \romannumeral\the\count@ pt\expandafter\endcsname
    \csname @\romannumeral\the\count@ pt\endcsname
  \csname #3\endcsname}%
\fi
\fi\endgroup

\begingroup\makeatletter\ifx\SetFigFont\undefined
\def\x#1#2#3#4#5#6#7\relax{\def\x{#1#2#3#4#5#6}}%
\expandafter\x\fmtname xxxxxx\relax \def\y{splain}%
\ifx\x\y   
\gdef\SetFigFont#1#2#3{%
  \ifnum #1<17\tiny\else \ifnum #1<20\small\else
  \ifnum #1<24\normalsize\else \ifnum #1<29\large\else
  \ifnum #1<34\Large\else \ifnum #1<41\LARGE\else
     \huge\fi\fi\fi\fi\fi\fi
  \csname #3\endcsname}%
\else
\gdef\SetFigFont#1#2#3{\begingroup
  \count@#1\relax \ifnum 25<\count@\count@25\fi
  \def\x{\endgroup\@setsize\SetFigFont{#2pt}}%
  \expandafter\x
    \csname \romannumeral\the\count@ pt\expandafter\endcsname
    \csname @\romannumeral\the\count@ pt\endcsname
  \csname #3\endcsname}%
\fi
\fi\endgroup

\newcommand{\Aiii}[3]{
\begin{picture}(2752,895)(1025,-4136)
\thinlines
\put(1201,-3961){\circle{336}}
\put(2401,-3961){\circle{336}}
\put(3601,-3961){\circle{336}}
\put(1369,-3961){\line( 1, 0){864}}
\put(2569,-3961){\line( 1, 0){864}}
\put(1141,-3436){\makebox(0,0)[lb]{\smash{\SetFigFont{6}{7.2}{rm}#1}}}
\put(2341,-3436){\makebox(0,0)[lb]{\smash{\SetFigFont{6}{7.2}{rm}#2}}}
\put(3541,-3436){\makebox(0,0)[lb]{\smash{\SetFigFont{6}{7.2}{rm}#3}}}
\end{picture}
}

\begingroup\makeatletter\ifx\SetFigFont\undefined
\def\x#1#2#3#4#5#6#7\relax{\def\x{#1#2#3#4#5#6}}%
\expandafter\x\fmtname xxxxxx\relax \def\y{splain}%
\ifx\x\y   
\gdef\SetFigFont#1#2#3{%
  \ifnum #1<17\tiny\else \ifnum #1<20\small\else
  \ifnum #1<24\normalsize\else \ifnum #1<29\large\else
  \ifnum #1<34\Large\else \ifnum #1<41\LARGE\else
     \huge\fi\fi\fi\fi\fi\fi
  \csname #3\endcsname}%
\else
\gdef\SetFigFont#1#2#3{\begingroup
  \count@#1\relax \ifnum 25<\count@\count@25\fi
  \def\x{\endgroup\@setsize\SetFigFont{#2pt}}%
  \expandafter\x
    \csname \romannumeral\the\count@ pt\expandafter\endcsname
    \csname @\romannumeral\the\count@ pt\endcsname
  \csname #3\endcsname}%
\fi
\fi\endgroup

\begingroup\makeatletter\ifx\SetFigFont\undefined
\def\x#1#2#3#4#5#6#7\relax{\def\x{#1#2#3#4#5#6}}%
\expandafter\x\fmtname xxxxxx\relax \def\y{splain}%
\ifx\x\y   
\gdef\SetFigFont#1#2#3{%
  \ifnum #1<17\tiny\else \ifnum #1<20\small\else
  \ifnum #1<24\normalsize\else \ifnum #1<29\large\else
  \ifnum #1<34\Large\else \ifnum #1<41\LARGE\else
     \huge\fi\fi\fi\fi\fi\fi
  \csname #3\endcsname}%
\else
\gdef\SetFigFont#1#2#3{\begingroup
  \count@#1\relax \ifnum 25<\count@\count@25\fi
  \def\x{\endgroup\@setsize\SetFigFont{#2pt}}%
  \expandafter\x
    \csname \romannumeral\the\count@ pt\expandafter\endcsname
    \csname @\romannumeral\the\count@ pt\endcsname
  \csname #3\endcsname}%
\fi
\fi\endgroup

\begingroup\makeatletter\ifx\SetFigFont\undefined
\def\x#1#2#3#4#5#6#7\relax{\def\x{#1#2#3#4#5#6}}%
\expandafter\x\fmtname xxxxxx\relax \def\y{splain}%
\ifx\x\y   
\gdef\SetFigFont#1#2#3{%
  \ifnum #1<17\tiny\else \ifnum #1<20\small\else
  \ifnum #1<24\normalsize\else \ifnum #1<29\large\else
  \ifnum #1<34\Large\else \ifnum #1<41\LARGE\else
     \huge\fi\fi\fi\fi\fi\fi
  \csname #3\endcsname}%
\else
\gdef\SetFigFont#1#2#3{\begingroup
  \count@#1\relax \ifnum 25<\count@\count@25\fi
  \def\x{\endgroup\@setsize\SetFigFont{#2pt}}%
  \expandafter\x
    \csname \romannumeral\the\count@ pt\expandafter\endcsname
    \csname @\romannumeral\the\count@ pt\endcsname
  \csname #3\endcsname}%
\fi
\fi\endgroup

\newcommand{\An}[3]{
\begin{picture}(3952,895)(1025,-4136)
\thinlines
\put(1201,-3961){\circle{336}}
\put(2401,-3961){\circle{336}}
\put(4771,-3961){\circle{336}}
\put(3436,-3961){\circle*{30}}
\put(3586,-3961){\circle*{30}}
\put(3736,-3961){\circle*{30}}
\put(1369,-3961){\line( 1, 0){870}}
\put(2569,-3961){\line( 1, 0){600}}
\put(4003,-3961){\line( 1, 0){600}}
\put(1141,-3436){\makebox(0,0)[lb]{\smash{\SetFigFont{6}{7.2}{rm}\mbox{$#1$}}}}
\put(2341,-3436){\makebox(0,0)[lb]{\smash{\SetFigFont{6}{7.2}{rm}\mbox{$#2$}}}}
\put(4741,-3436){\makebox(0,0)[lb]{\smash{\SetFigFont{6}{7.2}{rm}\mbox{$#3$}}}}
\end{picture}
}
\begingroup\makeatletter\ifx\SetFigFont\undefined
\def\x#1#2#3#4#5#6#7\relax{\def\x{#1#2#3#4#5#6}}%
\expandafter\x\fmtname xxxxxx\relax \def\y{splain}%
\ifx\x\y   
\gdef\SetFigFont#1#2#3{%
  \ifnum #1<17\tiny\else \ifnum #1<20\small\else
  \ifnum #1<24\normalsize\else \ifnum #1<29\large\else
  \ifnum #1<34\Large\else \ifnum #1<41\LARGE\else
     \huge\fi\fi\fi\fi\fi\fi
  \csname #3\endcsname}%
\else
\gdef\SetFigFont#1#2#3{\begingroup
  \count@#1\relax \ifnum 25<\count@\count@25\fi
  \def\x{\endgroup\@setsize\SetFigFont{#2pt}}%
  \expandafter\x
    \csname \romannumeral\the\count@ pt\expandafter\endcsname
    \csname @\romannumeral\the\count@ pt\endcsname
  \csname #3\endcsname}%
\fi
\fi\endgroup

\begingroup\makeatletter\ifx\SetFigFont\undefined
\def\x#1#2#3#4#5#6#7\relax{\def\x{#1#2#3#4#5#6}}%
\expandafter\x\fmtname xxxxxx\relax \def\y{splain}%
\ifx\x\y   
\gdef\SetFigFont#1#2#3{%
  \ifnum #1<17\tiny\else \ifnum #1<20\small\else
  \ifnum #1<24\normalsize\else \ifnum #1<29\large\else
  \ifnum #1<34\Large\else \ifnum #1<41\LARGE\else
     \huge\fi\fi\fi\fi\fi\fi
  \csname #3\endcsname}%
\else
\gdef\SetFigFont#1#2#3{\begingroup
  \count@#1\relax \ifnum 25<\count@\count@25\fi
  \def\x{\endgroup\@setsize\SetFigFont{#2pt}}%
  \expandafter\x
    \csname \romannumeral\the\count@ pt\expandafter\endcsname
    \csname @\romannumeral\the\count@ pt\endcsname
  \csname #3\endcsname}%
\fi
\fi\endgroup

\begingroup\makeatletter\ifx\SetFigFont\undefined
\def\x#1#2#3#4#5#6#7\relax{\def\x{#1#2#3#4#5#6}}%
\expandafter\x\fmtname xxxxxx\relax \def\y{splain}%
\ifx\x\y   
\gdef\SetFigFont#1#2#3{%
  \ifnum #1<17\tiny\else \ifnum #1<20\small\else
  \ifnum #1<24\normalsize\else \ifnum #1<29\large\else
  \ifnum #1<34\Large\else \ifnum #1<41\LARGE\else
     \huge\fi\fi\fi\fi\fi\fi
  \csname #3\endcsname}%
\else
\gdef\SetFigFont#1#2#3{\begingroup
  \count@#1\relax \ifnum 25<\count@\count@25\fi
  \def\x{\endgroup\@setsize\SetFigFont{#2pt}}%
  \expandafter\x
    \csname \romannumeral\the\count@ pt\expandafter\endcsname
    \csname @\romannumeral\the\count@ pt\endcsname
  \csname #3\endcsname}%
\fi
\fi\endgroup

\newcommand{\Ani}[4]{
\begin{picture}(5152,895)(1025,-4136)
\thinlines
\put(1201,-3961){\circle{336}}
\put(2401,-3961){\circle{336}}
\put(4771,-3961){\circle{336}}
\put(5971,-3961){\circle{336}}
\put(3436,-3961){\circle*{30}}
\put(3586,-3961){\circle*{30}}
\put(3736,-3961){\circle*{30}}
\put(1369,-3961){\line( 1, 0){870}}
\put(2569,-3961){\line( 1, 0){600}}
\put(4003,-3961){\line( 1, 0){600}}
\put(4969,-3961){\line( 1, 0){870}}
\put(1141,-3436){\makebox(0,0)[lb]{\smash{\SetFigFont{6}{7.2}{rm}\mbox{$#1$}}}}
\put(2341,-3436){\makebox(0,0)[lb]{\smash{\SetFigFont{6}{7.2}{rm}\mbox{$#2$}}}}
\put(4741,-3436){\makebox(0,0)[lb]{\smash{\SetFigFont{6}{7.2}{rm}\mbox{$#3$}}}}
\put(5941,-3436){\makebox(0,0)[lb]{\smash{\SetFigFont{6}{7.2}{rm}\mbox{$#4$}}}}
\end{picture}
}
\begingroup\makeatletter\ifx\SetFigFont\undefined
\def\x#1#2#3#4#5#6#7\relax{\def\x{#1#2#3#4#5#6}}%
\expandafter\x\fmtname xxxxxx\relax \def\y{splain}%
\ifx\x\y   
\gdef\SetFigFont#1#2#3{%
  \ifnum #1<17\tiny\else \ifnum #1<20\small\else
  \ifnum #1<24\normalsize\else \ifnum #1<29\large\else
  \ifnum #1<34\Large\else \ifnum #1<41\LARGE\else
     \huge\fi\fi\fi\fi\fi\fi
  \csname #3\endcsname}%
\else
\gdef\SetFigFont#1#2#3{\begingroup
  \count@#1\relax \ifnum 25<\count@\count@25\fi
  \def\x{\endgroup\@setsize\SetFigFont{#2pt}}%
  \expandafter\x
    \csname \romannumeral\the\count@ pt\expandafter\endcsname
    \csname @\romannumeral\the\count@ pt\endcsname
  \csname #3\endcsname}%
\fi
\fi\endgroup

\newcommand{\Anii}[7]{
\begin{picture}(9952,895)(1025,-4136)
\thinlines
\put(1201,-3961){\circle{336}}
\put(2401,-3961){\circle{336}}
\put(3436,-3961){\circle*{30}}
\put(3586,-3961){\circle*{30}}
\put(3736,-3961){\circle*{30}}
\put(4771,-3961){\circle{336}}
\put(5971,-3961){\circle{336}}
\put(7171,-3961){\circle{336}}
\put(8206,-3961){\circle*{30}}
\put(8356,-3961){\circle*{30}}
\put(8506,-3961){\circle*{30}}
\put(9571,-3961){\circle{336}}
\put(10771,-3961){\circle{336}}

\put(1369,-3961){\line( 1, 0){864}}
\put(2569,-3961){\line( 1, 0){600}}
\put(4003,-3961){\line( 1, 0){600}}
\put(4969,-3961){\line( 1, 0){864}}
\put(6169,-3961){\line( 1, 0){864}}
\put(7369,-3961){\line( 1, 0){600}}
\put(8803,-3961){\line( 1, 0){600}}
\put(9769,-3961){\line( 1, 0){864}}

\put(1141,-3436){\makebox(0,0)[lb]{\smash{\SetFigFont{6}{7.2}{rm}\mbox{$#1$}}}}
\put(2341,-3436){\makebox(0,0)[lb]{\smash{\SetFigFont{6}{7.2}{rm}\mbox{$#2$}}}}
\put(4741,-3436){\makebox(0,0)[lb]{\smash{\SetFigFont{6}{7.2}{rm}\mbox{$#3$}}}}
\put(5941,-3436){\makebox(0,0)[lb]{\smash{\SetFigFont{6}{7.2}{rm}\mbox{$#4$}}}}
\put(7141,-3436){\makebox(0,0)[lb]{\smash{\SetFigFont{6}{7.2}{rm}\mbox{$#5$}}}}
\put(9541,-3436){\makebox(0,0)[lb]{\smash{\SetFigFont{6}{7.2}{rm}\mbox{$#6$}}}}
\put(10741,-3436){\makebox(0,0)[lb]{\smash{\SetFigFont{6}{7.2}{rm}\mbox{$#7$}}}}
\end{picture}
}
\begingroup\makeatletter\ifx\SetFigFont\undefined
\def\x#1#2#3#4#5#6#7\relax{\def\x{#1#2#3#4#5#6}}%
\expandafter\x\fmtname xxxxxx\relax \def\y{splain}%
\ifx\x\y   
\gdef\SetFigFont#1#2#3{%
  \ifnum #1<17\tiny\else \ifnum #1<20\small\else
  \ifnum #1<24\normalsize\else \ifnum #1<29\large\else
  \ifnum #1<34\Large\else \ifnum #1<41\LARGE\else
     \huge\fi\fi\fi\fi\fi\fi
  \csname #3\endcsname}%
\else
\gdef\SetFigFont#1#2#3{\begingroup
  \count@#1\relax \ifnum 25<\count@\count@25\fi
  \def\x{\endgroup\@setsize\SetFigFont{#2pt}}%
  \expandafter\x
    \csname \romannumeral\the\count@ pt\expandafter\endcsname
    \csname @\romannumeral\the\count@ pt\endcsname
  \csname #3\endcsname}%
\fi
\fi\endgroup

\begingroup\makeatletter\ifx\SetFigFont\undefined
\def\x#1#2#3#4#5#6#7\relax{\def\x{#1#2#3#4#5#6}}%
\expandafter\x\fmtname xxxxxx\relax \def\y{splain}%
\ifx\x\y   
\gdef\SetFigFont#1#2#3{%
  \ifnum #1<17\tiny\else \ifnum #1<20\small\else
  \ifnum #1<24\normalsize\else \ifnum #1<29\large\else
  \ifnum #1<34\Large\else \ifnum #1<41\LARGE\else
     \huge\fi\fi\fi\fi\fi\fi
  \csname #3\endcsname}%
\else
\gdef\SetFigFont#1#2#3{\begingroup
  \count@#1\relax \ifnum 25<\count@\count@25\fi
  \def\x{\endgroup\@setsize\SetFigFont{#2pt}}%
  \expandafter\x
    \csname \romannumeral\the\count@ pt\expandafter\endcsname
    \csname @\romannumeral\the\count@ pt\endcsname
  \csname #3\endcsname}%
\fi
\fi\endgroup

\begingroup\makeatletter\ifx\SetFigFont\undefined
\def\x#1#2#3#4#5#6#7\relax{\def\x{#1#2#3#4#5#6}}%
\expandafter\x\fmtname xxxxxx\relax \def\y{splain}%
\ifx\x\y   
\gdef\SetFigFont#1#2#3{%
  \ifnum #1<17\tiny\else \ifnum #1<20\small\else
  \ifnum #1<24\normalsize\else \ifnum #1<29\large\else
  \ifnum #1<34\Large\else \ifnum #1<41\LARGE\else
     \huge\fi\fi\fi\fi\fi\fi
  \csname #3\endcsname}%
\else
\gdef\SetFigFont#1#2#3{\begingroup
  \count@#1\relax \ifnum 25<\count@\count@25\fi
  \def\x{\endgroup\@setsize\SetFigFont{#2pt}}%
  \expandafter\x
    \csname \romannumeral\the\count@ pt\expandafter\endcsname
    \csname @\romannumeral\the\count@ pt\endcsname
  \csname #3\endcsname}%
\fi
\fi\endgroup

\begingroup\makeatletter\ifx\SetFigFont\undefined
\def\x#1#2#3#4#5#6#7\relax{\def\x{#1#2#3#4#5#6}}%
\expandafter\x\fmtname xxxxxx\relax \def\y{splain}%
\ifx\x\y   
\gdef\SetFigFont#1#2#3{%
  \ifnum #1<17\tiny\else \ifnum #1<20\small\else
  \ifnum #1<24\normalsize\else \ifnum #1<29\large\else
  \ifnum #1<34\Large\else \ifnum #1<41\LARGE\else
     \huge\fi\fi\fi\fi\fi\fi
  \csname #3\endcsname}%
\else
\gdef\SetFigFont#1#2#3{\begingroup
  \count@#1\relax \ifnum 25<\count@\count@25\fi
  \def\x{\endgroup\@setsize\SetFigFont{#2pt}}%
  \expandafter\x
    \csname \romannumeral\the\count@ pt\expandafter\endcsname
    \csname @\romannumeral\the\count@ pt\endcsname
  \csname #3\endcsname}%
\fi
\fi\endgroup

\begingroup\makeatletter\ifx\SetFigFont\undefined
\def\x#1#2#3#4#5#6#7\relax{\def\x{#1#2#3#4#5#6}}%
\expandafter\x\fmtname xxxxxx\relax \def\y{splain}%
\ifx\x\y   
\gdef\SetFigFont#1#2#3{%
  \ifnum #1<17\tiny\else \ifnum #1<20\small\else
  \ifnum #1<24\normalsize\else \ifnum #1<29\large\else
  \ifnum #1<34\Large\else \ifnum #1<41\LARGE\else
     \huge\fi\fi\fi\fi\fi\fi
  \csname #3\endcsname}%
\else
\gdef\SetFigFont#1#2#3{\begingroup
  \count@#1\relax \ifnum 25<\count@\count@25\fi
  \def\x{\endgroup\@setsize\SetFigFont{#2pt}}%
  \expandafter\x
    \csname \romannumeral\the\count@ pt\expandafter\endcsname
    \csname @\romannumeral\the\count@ pt\endcsname
  \csname #3\endcsname}%
\fi
\fi\endgroup

\begingroup\makeatletter\ifx\SetFigFont\undefined
\def\x#1#2#3#4#5#6#7\relax{\def\x{#1#2#3#4#5#6}}%
\expandafter\x\fmtname xxxxxx\relax \def\y{splain}%
\ifx\x\y   
\gdef\SetFigFont#1#2#3{%
  \ifnum #1<17\tiny\else \ifnum #1<20\small\else
  \ifnum #1<24\normalsize\else \ifnum #1<29\large\else
  \ifnum #1<34\Large\else \ifnum #1<41\LARGE\else
     \huge\fi\fi\fi\fi\fi\fi
  \csname #3\endcsname}%
\else
\gdef\SetFigFont#1#2#3{\begingroup
  \count@#1\relax \ifnum 25<\count@\count@25\fi
  \def\x{\endgroup\@setsize\SetFigFont{#2pt}}%
  \expandafter\x
    \csname \romannumeral\the\count@ pt\expandafter\endcsname
    \csname @\romannumeral\the\count@ pt\endcsname
  \csname #3\endcsname}%
\fi
\fi\endgroup

\begingroup\makeatletter\ifx\SetFigFont\undefined
\def\x#1#2#3#4#5#6#7\relax{\def\x{#1#2#3#4#5#6}}%
\expandafter\x\fmtname xxxxxx\relax \def\y{splain}%
\ifx\x\y   
\gdef\SetFigFont#1#2#3{%
  \ifnum #1<17\tiny\else \ifnum #1<20\small\else
  \ifnum #1<24\normalsize\else \ifnum #1<29\large\else
  \ifnum #1<34\Large\else \ifnum #1<41\LARGE\else
     \huge\fi\fi\fi\fi\fi\fi
  \csname #3\endcsname}%
\else
\gdef\SetFigFont#1#2#3{\begingroup
  \count@#1\relax \ifnum 25<\count@\count@25\fi
  \def\x{\endgroup\@setsize\SetFigFont{#2pt}}%
  \expandafter\x
    \csname \romannumeral\the\count@ pt\expandafter\endcsname
    \csname @\romannumeral\the\count@ pt\endcsname
  \csname #3\endcsname}%
\fi
\fi\endgroup

\newcommand{\Biii}[3]{
\begin{picture}(2752,895)(7025,-4136)
\thinlines
\put(7201,-3961){\circle{336}}
\put(8401,-3961){\circle{336}}
\put(9601,-3961){\circle*{336}}
\put(7369,-3961){\line( 1, 0){864}}
\put(8569,-3911){\line( 1, 0){864}}
\put(8569,-4011){\line( 1, 0){864}}
\put(7141,-3436){\makebox(0,0)[lb]{\smash{\SetFigFont{6}{7.2}{rm}#1}}}
\put(8341,-3436){\makebox(0,0)[lb]{\smash{\SetFigFont{6}{7.2}{rm}#2}}}
\put(9541,-3436){\makebox(0,0)[lb]{\smash{\SetFigFont{6}{7.2}{rm}#3}}}
\end{picture}
}
\begingroup\makeatletter\ifx\SetFigFont\undefined
\def\x#1#2#3#4#5#6#7\relax{\def\x{#1#2#3#4#5#6}}%
\expandafter\x\fmtname xxxxxx\relax \def\y{splain}%
\ifx\x\y   
\gdef\SetFigFont#1#2#3{%
  \ifnum #1<17\tiny\else \ifnum #1<20\small\else
  \ifnum #1<24\normalsize\else \ifnum #1<29\large\else
  \ifnum #1<34\Large\else \ifnum #1<41\LARGE\else
     \huge\fi\fi\fi\fi\fi\fi
  \csname #3\endcsname}%
\else
\gdef\SetFigFont#1#2#3{\begingroup
  \count@#1\relax \ifnum 25<\count@\count@25\fi
  \def\x{\endgroup\@setsize\SetFigFont{#2pt}}%
  \expandafter\x
    \csname \romannumeral\the\count@ pt\expandafter\endcsname
    \csname @\romannumeral\the\count@ pt\endcsname
  \csname #3\endcsname}%
\fi
\fi\endgroup

\begingroup\makeatletter\ifx\SetFigFont\undefined
\def\x#1#2#3#4#5#6#7\relax{\def\x{#1#2#3#4#5#6}}%
\expandafter\x\fmtname xxxxxx\relax \def\y{splain}%
\ifx\x\y   
\gdef\SetFigFont#1#2#3{%
  \ifnum #1<17\tiny\else \ifnum #1<20\small\else
  \ifnum #1<24\normalsize\else \ifnum #1<29\large\else
  \ifnum #1<34\Large\else \ifnum #1<41\LARGE\else
     \huge\fi\fi\fi\fi\fi\fi
  \csname #3\endcsname}%
\else
\gdef\SetFigFont#1#2#3{\begingroup
  \count@#1\relax \ifnum 25<\count@\count@25\fi
  \def\x{\endgroup\@setsize\SetFigFont{#2pt}}%
  \expandafter\x
    \csname \romannumeral\the\count@ pt\expandafter\endcsname
    \csname @\romannumeral\the\count@ pt\endcsname
  \csname #3\endcsname}%
\fi
\fi\endgroup

\newcommand{\Bn}[4]{
\begin{picture}(5152,895)(1025,-4136)
\thinlines
\put(1201,-3961){\circle{336}}
\put(2401,-3961){\circle{336}}
\put(4801,-3961){\circle{336}}
\put(6001,-3961){\circle*{336}}
\put(3436,-3946){\circle*{30}}
\put(3586,-3946){\circle*{30}}
\put(3736,-3946){\circle*{30}}
\put(1369,-3961){\line( 1, 0){864}}
\put(2569,-3961){\line( 1, 0){600}}
\put(4003,-3961){\line( 1, 0){600}}
\put(4969,-3911){\line( 1, 0){864}}
\put(4969,-4011){\line( 1, 0){864}}
\put(1141,-3436){\makebox(0,0)[lb]{\smash{\SetFigFont{6}{7.2}{rm}\mbox{$#1$}}}}
\put(2341,-3436){\makebox(0,0)[lb]{\smash{\SetFigFont{6}{7.2}{rm}\mbox{$#2$}}}}
\put(4741,-3436){\makebox(0,0)[lb]{\smash{\SetFigFont{6}{7.2}{rm}\mbox{$#3$}}}}
\put(5941,-3436){\makebox(0,0)[lb]{\smash{\SetFigFont{6}{7.2}{rm}\mbox{$#4$}}}}
\end{picture}
}
\begingroup\makeatletter\ifx\SetFigFont\undefined
\def\x#1#2#3#4#5#6#7\relax{\def\x{#1#2#3#4#5#6}}%
\expandafter\x\fmtname xxxxxx\relax \def\y{splain}%
\ifx\x\y   
\gdef\SetFigFont#1#2#3{%
  \ifnum #1<17\tiny\else \ifnum #1<20\small\else
  \ifnum #1<24\normalsize\else \ifnum #1<29\large\else
  \ifnum #1<34\Large\else \ifnum #1<41\LARGE\else
     \huge\fi\fi\fi\fi\fi\fi
  \csname #3\endcsname}%
\else
\gdef\SetFigFont#1#2#3{\begingroup
  \count@#1\relax \ifnum 25<\count@\count@25\fi
  \def\x{\endgroup\@setsize\SetFigFont{#2pt}}%
  \expandafter\x
    \csname \romannumeral\the\count@ pt\expandafter\endcsname
    \csname @\romannumeral\the\count@ pt\endcsname
  \csname #3\endcsname}%
\fi
\fi\endgroup

\begingroup\makeatletter\ifx\SetFigFont\undefined
\def\x#1#2#3#4#5#6#7\relax{\def\x{#1#2#3#4#5#6}}%
\expandafter\x\fmtname xxxxxx\relax \def\y{splain}%
\ifx\x\y   
\gdef\SetFigFont#1#2#3{%
  \ifnum #1<17\tiny\else \ifnum #1<20\small\else
  \ifnum #1<24\normalsize\else \ifnum #1<29\large\else
  \ifnum #1<34\Large\else \ifnum #1<41\LARGE\else
     \huge\fi\fi\fi\fi\fi\fi
  \csname #3\endcsname}%
\else
\gdef\SetFigFont#1#2#3{\begingroup
  \count@#1\relax \ifnum 25<\count@\count@25\fi
  \def\x{\endgroup\@setsize\SetFigFont{#2pt}}%
  \expandafter\x
    \csname \romannumeral\the\count@ pt\expandafter\endcsname
    \csname @\romannumeral\the\count@ pt\endcsname
  \csname #3\endcsname}%
\fi
\fi\endgroup

\begingroup\makeatletter\ifx\SetFigFont\undefined
\def\x#1#2#3#4#5#6#7\relax{\def\x{#1#2#3#4#5#6}}%
\expandafter\x\fmtname xxxxxx\relax \def\y{splain}%
\ifx\x\y   
\gdef\SetFigFont#1#2#3{%
  \ifnum #1<17\tiny\else \ifnum #1<20\small\else
  \ifnum #1<24\normalsize\else \ifnum #1<29\large\else
  \ifnum #1<34\Large\else \ifnum #1<41\LARGE\else
     \huge\fi\fi\fi\fi\fi\fi
  \csname #3\endcsname}%
\else
\gdef\SetFigFont#1#2#3{\begingroup
  \count@#1\relax \ifnum 25<\count@\count@25\fi
  \def\x{\endgroup\@setsize\SetFigFont{#2pt}}%
  \expandafter\x
    \csname \romannumeral\the\count@ pt\expandafter\endcsname
    \csname @\romannumeral\the\count@ pt\endcsname
  \csname #3\endcsname}%
\fi
\fi\endgroup

\newcommand{\Bnii}[6]{
\begin{picture}(7552,895)(1025,-4136)
\thinlines
\put(1201,-3961){\circle{336}}
\put(2401,-3961){\circle{336}}
\put(3601,-3961){\circle{336}}
\put(4801,-3961){\circle{336}}
\put(7201,-3961){\circle{336}}
\put(8401,-3961){\circle*{336}}

\put(1369,-3961){\line( 1, 0){864}}
\put(2569,-3961){\line( 1, 0){864}}
\put(3769,-3961){\line( 1, 0){864}}

\put(4903,-3961){\line( 1, 0){600}}
\put(6403,-3961){\line( 1, 0){600}}
\put(7369,-3911){\line( 1, 0){864}}
\put(7369,-4011){\line( 1, 0){864}}

\put(5836,-3946){\circle*{30}}
\put(5986,-3946){\circle*{30}}
\put(6136,-3946){\circle*{30}}

\put(1141,-3436){\makebox(0,0)[lb]{\smash{\SetFigFont{6}{7.2}{rm}\mbox{$#1$}}}}
\put(2341,-3436){\makebox(0,0)[lb]{\smash{\SetFigFont{6}{7.2}{rm}\mbox{$#2$}}}}
\put(3541,-3436){\makebox(0,0)[lb]{\smash{\SetFigFont{6}{7.2}{rm}\mbox{$#3$}}}}
\put(4741,-3436){\makebox(0,0)[lb]{\smash{\SetFigFont{6}{7.2}{rm}\mbox{$#4$}}}}
\put(7141,-3436){\makebox(0,0)[lb]{\smash{\SetFigFont{6}{7.2}{rm}\mbox{$#5$}}}}
\put(8341,-3436){\makebox(0,0)[lb]{\smash{\SetFigFont{6}{7.2}{rm}\mbox{$#6$}}}}

\end{picture}
}
\begingroup\makeatletter\ifx\SetFigFont\undefined
\def\x#1#2#3#4#5#6#7\relax{\def\x{#1#2#3#4#5#6}}%
\expandafter\x\fmtname xxxxxx\relax \def\y{splain}%
\ifx\x\y   
\gdef\SetFigFont#1#2#3{%
  \ifnum #1<17\tiny\else \ifnum #1<20\small\else
  \ifnum #1<24\normalsize\else \ifnum #1<29\large\else
  \ifnum #1<34\Large\else \ifnum #1<41\LARGE\else
     \huge\fi\fi\fi\fi\fi\fi
  \csname #3\endcsname}%
\else
\gdef\SetFigFont#1#2#3{\begingroup
  \count@#1\relax \ifnum 25<\count@\count@25\fi
  \def\x{\endgroup\@setsize\SetFigFont{#2pt}}%
  \expandafter\x
    \csname \romannumeral\the\count@ pt\expandafter\endcsname
    \csname @\romannumeral\the\count@ pt\endcsname
  \csname #3\endcsname}%
\fi
\fi\endgroup

\begingroup\makeatletter\ifx\SetFigFont\undefined
\def\x#1#2#3#4#5#6#7\relax{\def\x{#1#2#3#4#5#6}}%
\expandafter\x\fmtname xxxxxx\relax \def\y{splain}%
\ifx\x\y   
\gdef\SetFigFont#1#2#3{%
  \ifnum #1<17\tiny\else \ifnum #1<20\small\else
  \ifnum #1<24\normalsize\else \ifnum #1<29\large\else
  \ifnum #1<34\Large\else \ifnum #1<41\LARGE\else
     \huge\fi\fi\fi\fi\fi\fi
  \csname #3\endcsname}%
\else
\gdef\SetFigFont#1#2#3{\begingroup
  \count@#1\relax \ifnum 25<\count@\count@25\fi
  \def\x{\endgroup\@setsize\SetFigFont{#2pt}}%
  \expandafter\x
    \csname \romannumeral\the\count@ pt\expandafter\endcsname
    \csname @\romannumeral\the\count@ pt\endcsname
  \csname #3\endcsname}%
\fi
\fi\endgroup

\begingroup\makeatletter\ifx\SetFigFont\undefined
\def\x#1#2#3#4#5#6#7\relax{\def\x{#1#2#3#4#5#6}}%
\expandafter\x\fmtname xxxxxx\relax \def\y{splain}%
\ifx\x\y   
\gdef\SetFigFont#1#2#3{%
  \ifnum #1<17\tiny\else \ifnum #1<20\small\else
  \ifnum #1<24\normalsize\else \ifnum #1<29\large\else
  \ifnum #1<34\Large\else \ifnum #1<41\LARGE\else
     \huge\fi\fi\fi\fi\fi\fi
  \csname #3\endcsname}%
\else
\gdef\SetFigFont#1#2#3{\begingroup
  \count@#1\relax \ifnum 25<\count@\count@25\fi
  \def\x{\endgroup\@setsize\SetFigFont{#2pt}}%
  \expandafter\x
    \csname \romannumeral\the\count@ pt\expandafter\endcsname
    \csname @\romannumeral\the\count@ pt\endcsname
  \csname #3\endcsname}%
\fi
\fi\endgroup

\begingroup\makeatletter\ifx\SetFigFont\undefined
\def\x#1#2#3#4#5#6#7\relax{\def\x{#1#2#3#4#5#6}}%
\expandafter\x\fmtname xxxxxx\relax \def\y{splain}%
\ifx\x\y   
\gdef\SetFigFont#1#2#3{%
  \ifnum #1<17\tiny\else \ifnum #1<20\small\else
  \ifnum #1<24\normalsize\else \ifnum #1<29\large\else
  \ifnum #1<34\Large\else \ifnum #1<41\LARGE\else
     \huge\fi\fi\fi\fi\fi\fi
  \csname #3\endcsname}%
\else
\gdef\SetFigFont#1#2#3{\begingroup
  \count@#1\relax \ifnum 25<\count@\count@25\fi
  \def\x{\endgroup\@setsize\SetFigFont{#2pt}}%
  \expandafter\x
    \csname \romannumeral\the\count@ pt\expandafter\endcsname
    \csname @\romannumeral\the\count@ pt\endcsname
  \csname #3\endcsname}%
\fi
\fi\endgroup

\begingroup\makeatletter\ifx\SetFigFont\undefined
\def\x#1#2#3#4#5#6#7\relax{\def\x{#1#2#3#4#5#6}}%
\expandafter\x\fmtname xxxxxx\relax \def\y{splain}%
\ifx\x\y   
\gdef\SetFigFont#1#2#3{%
  \ifnum #1<17\tiny\else \ifnum #1<20\small\else
  \ifnum #1<24\normalsize\else \ifnum #1<29\large\else
  \ifnum #1<34\Large\else \ifnum #1<41\LARGE\else
     \huge\fi\fi\fi\fi\fi\fi
  \csname #3\endcsname}%
\else
\gdef\SetFigFont#1#2#3{\begingroup
  \count@#1\relax \ifnum 25<\count@\count@25\fi
  \def\x{\endgroup\@setsize\SetFigFont{#2pt}}%
  \expandafter\x
    \csname \romannumeral\the\count@ pt\expandafter\endcsname
    \csname @\romannumeral\the\count@ pt\endcsname
  \csname #3\endcsname}%
\fi
\fi\endgroup

\begingroup\makeatletter\ifx\SetFigFont\undefined
\def\x#1#2#3#4#5#6#7\relax{\def\x{#1#2#3#4#5#6}}%
\expandafter\x\fmtname xxxxxx\relax \def\y{splain}%
\ifx\x\y   
\gdef\SetFigFont#1#2#3{%
  \ifnum #1<17\tiny\else \ifnum #1<20\small\else
  \ifnum #1<24\normalsize\else \ifnum #1<29\large\else
  \ifnum #1<34\Large\else \ifnum #1<41\LARGE\else
     \huge\fi\fi\fi\fi\fi\fi
  \csname #3\endcsname}%
\else
\gdef\SetFigFont#1#2#3{\begingroup
  \count@#1\relax \ifnum 25<\count@\count@25\fi
  \def\x{\endgroup\@setsize\SetFigFont{#2pt}}%
  \expandafter\x
    \csname \romannumeral\the\count@ pt\expandafter\endcsname
    \csname @\romannumeral\the\count@ pt\endcsname
  \csname #3\endcsname}%
\fi
\fi\endgroup

\begingroup\makeatletter\ifx\SetFigFont\undefined
\def\x#1#2#3#4#5#6#7\relax{\def\x{#1#2#3#4#5#6}}%
\expandafter\x\fmtname xxxxxx\relax \def\y{splain}%
\ifx\x\y   
\gdef\SetFigFont#1#2#3{%
  \ifnum #1<17\tiny\else \ifnum #1<20\small\else
  \ifnum #1<24\normalsize\else \ifnum #1<29\large\else
  \ifnum #1<34\Large\else \ifnum #1<41\LARGE\else
     \huge\fi\fi\fi\fi\fi\fi
  \csname #3\endcsname}%
\else
\gdef\SetFigFont#1#2#3{\begingroup
  \count@#1\relax \ifnum 25<\count@\count@25\fi
  \def\x{\endgroup\@setsize\SetFigFont{#2pt}}%
  \expandafter\x
    \csname \romannumeral\the\count@ pt\expandafter\endcsname
    \csname @\romannumeral\the\count@ pt\endcsname
  \csname #3\endcsname}%
\fi
\fi\endgroup

\begingroup\makeatletter\ifx\SetFigFont\undefined
\def\x#1#2#3#4#5#6#7\relax{\def\x{#1#2#3#4#5#6}}%
\expandafter\x\fmtname xxxxxx\relax \def\y{splain}%
\ifx\x\y   
\gdef\SetFigFont#1#2#3{%
  \ifnum #1<17\tiny\else \ifnum #1<20\small\else
  \ifnum #1<24\normalsize\else \ifnum #1<29\large\else
  \ifnum #1<34\Large\else \ifnum #1<41\LARGE\else
     \huge\fi\fi\fi\fi\fi\fi
  \csname #3\endcsname}%
\else
\gdef\SetFigFont#1#2#3{\begingroup
  \count@#1\relax \ifnum 25<\count@\count@25\fi
  \def\x{\endgroup\@setsize\SetFigFont{#2pt}}%
  \expandafter\x
    \csname \romannumeral\the\count@ pt\expandafter\endcsname
    \csname @\romannumeral\the\count@ pt\endcsname
  \csname #3\endcsname}%
\fi
\fi\endgroup

\begingroup\makeatletter\ifx\SetFigFont\undefined
\def\x#1#2#3#4#5#6#7\relax{\def\x{#1#2#3#4#5#6}}%
\expandafter\x\fmtname xxxxxx\relax \def\y{splain}%
\ifx\x\y   
\gdef\SetFigFont#1#2#3{%
  \ifnum #1<17\tiny\else \ifnum #1<20\small\else
  \ifnum #1<24\normalsize\else \ifnum #1<29\large\else
  \ifnum #1<34\Large\else \ifnum #1<41\LARGE\else
     \huge\fi\fi\fi\fi\fi\fi
  \csname #3\endcsname}%
\else
\gdef\SetFigFont#1#2#3{\begingroup
  \count@#1\relax \ifnum 25<\count@\count@25\fi
  \def\x{\endgroup\@setsize\SetFigFont{#2pt}}%
  \expandafter\x
    \csname \romannumeral\the\count@ pt\expandafter\endcsname
    \csname @\romannumeral\the\count@ pt\endcsname
  \csname #3\endcsname}%
\fi
\fi\endgroup

\begingroup\makeatletter\ifx\SetFigFont\undefined
\def\x#1#2#3#4#5#6#7\relax{\def\x{#1#2#3#4#5#6}}%
\expandafter\x\fmtname xxxxxx\relax \def\y{splain}%
\ifx\x\y   
\gdef\SetFigFont#1#2#3{%
  \ifnum #1<17\tiny\else \ifnum #1<20\small\else
  \ifnum #1<24\normalsize\else \ifnum #1<29\large\else
  \ifnum #1<34\Large\else \ifnum #1<41\LARGE\else
     \huge\fi\fi\fi\fi\fi\fi
  \csname #3\endcsname}%
\else
\gdef\SetFigFont#1#2#3{\begingroup
  \count@#1\relax \ifnum 25<\count@\count@25\fi
  \def\x{\endgroup\@setsize\SetFigFont{#2pt}}%
  \expandafter\x
    \csname \romannumeral\the\count@ pt\expandafter\endcsname
    \csname @\romannumeral\the\count@ pt\endcsname
  \csname #3\endcsname}%
\fi
\fi\endgroup

\newcommand{\Cn}[4]{
\begin{picture}(5152,895)(1025,-4136)
\thinlines
\put(1201,-3961){\circle*{336}}
\put(2401,-3961){\circle*{336}}
\put(4801,-3961){\circle*{336}}
\put(6001,-3961){\circle{336}}
\put(3436,-3946){\circle*{30}}
\put(3586,-3946){\circle*{30}}
\put(3736,-3946){\circle*{30}}
\put(1369,-3961){\line( 1, 0){864}}
\put(2569,-3961){\line( 1, 0){600}}
\put(4003,-3961){\line( 1, 0){600}}
\put(4969,-3911){\line( 1, 0){864}}
\put(4969,-4011){\line( 1, 0){864}}
\put(1141,-3436){\makebox(0,0)[lb]{\smash{\SetFigFont{6}{7.2}{rm}\mbox{$#1$}}}}
\put(2341,-3436){\makebox(0,0)[lb]{\smash{\SetFigFont{6}{7.2}{rm}\mbox{$#2$}}}}
\put(4741,-3436){\makebox(0,0)[lb]{\smash{\SetFigFont{6}{7.2}{rm}\mbox{$#3$}}}}
\put(5941,-3436){\makebox(0,0)[lb]{\smash{\SetFigFont{6}{7.2}{rm}\mbox{$#4$}}}}
\end{picture}
}
\begingroup\makeatletter\ifx\SetFigFont\undefined
\def\x#1#2#3#4#5#6#7\relax{\def\x{#1#2#3#4#5#6}}%
\expandafter\x\fmtname xxxxxx\relax \def\y{splain}%
\ifx\x\y   
\gdef\SetFigFont#1#2#3{%
  \ifnum #1<17\tiny\else \ifnum #1<20\small\else
  \ifnum #1<24\normalsize\else \ifnum #1<29\large\else
  \ifnum #1<34\Large\else \ifnum #1<41\LARGE\else
     \huge\fi\fi\fi\fi\fi\fi
  \csname #3\endcsname}%
\else
\gdef\SetFigFont#1#2#3{\begingroup
  \count@#1\relax \ifnum 25<\count@\count@25\fi
  \def\x{\endgroup\@setsize\SetFigFont{#2pt}}%
  \expandafter\x
    \csname \romannumeral\the\count@ pt\expandafter\endcsname
    \csname @\romannumeral\the\count@ pt\endcsname
  \csname #3\endcsname}%
\fi
\fi\endgroup

\begingroup\makeatletter\ifx\SetFigFont\undefined
\def\x#1#2#3#4#5#6#7\relax{\def\x{#1#2#3#4#5#6}}%
\expandafter\x\fmtname xxxxxx\relax \def\y{splain}%
\ifx\x\y   
\gdef\SetFigFont#1#2#3{%
  \ifnum #1<17\tiny\else \ifnum #1<20\small\else
  \ifnum #1<24\normalsize\else \ifnum #1<29\large\else
  \ifnum #1<34\Large\else \ifnum #1<41\LARGE\else
     \huge\fi\fi\fi\fi\fi\fi
  \csname #3\endcsname}%
\else
\gdef\SetFigFont#1#2#3{\begingroup
  \count@#1\relax \ifnum 25<\count@\count@25\fi
  \def\x{\endgroup\@setsize\SetFigFont{#2pt}}%
  \expandafter\x
    \csname \romannumeral\the\count@ pt\expandafter\endcsname
    \csname @\romannumeral\the\count@ pt\endcsname
  \csname #3\endcsname}%
\fi
\fi\endgroup

\begingroup\makeatletter\ifx\SetFigFont\undefined
\def\x#1#2#3#4#5#6#7\relax{\def\x{#1#2#3#4#5#6}}%
\expandafter\x\fmtname xxxxxx\relax \def\y{splain}%
\ifx\x\y   
\gdef\SetFigFont#1#2#3{%
  \ifnum #1<17\tiny\else \ifnum #1<20\small\else
  \ifnum #1<24\normalsize\else \ifnum #1<29\large\else
  \ifnum #1<34\Large\else \ifnum #1<41\LARGE\else
     \huge\fi\fi\fi\fi\fi\fi
  \csname #3\endcsname}%
\else
\gdef\SetFigFont#1#2#3{\begingroup
  \count@#1\relax \ifnum 25<\count@\count@25\fi
  \def\x{\endgroup\@setsize\SetFigFont{#2pt}}%
  \expandafter\x
    \csname \romannumeral\the\count@ pt\expandafter\endcsname
    \csname @\romannumeral\the\count@ pt\endcsname
  \csname #3\endcsname}%
\fi
\fi\endgroup

\newcommand{\Cnii}[6]{
\begin{picture}(7552,895)(1025,-4136)
\thinlines
\put(1201,-3961){\circle*{336}}
\put(2401,-3961){\circle*{336}}
\put(3601,-3961){\circle*{336}}
\put(4801,-3961){\circle*{336}}
\put(7201,-3961){\circle*{336}}
\put(8401,-3961){\circle{336}}

\put(1369,-3961){\line( 1, 0){864}}
\put(2569,-3961){\line( 1, 0){864}}
\put(3769,-3961){\line( 1, 0){864}}

\put(4903,-3961){\line( 1, 0){600}}
\put(6403,-3961){\line( 1, 0){600}}
\put(7369,-3911){\line( 1, 0){864}}
\put(7369,-4011){\line( 1, 0){864}}

\put(5836,-3946){\circle*{30}}
\put(5986,-3946){\circle*{30}}
\put(6136,-3946){\circle*{30}}

\put(1141,-3436){\makebox(0,0)[lb]{\smash{\SetFigFont{6}{7.2}{rm}#1}}}
\put(2341,-3436){\makebox(0,0)[lb]{\smash{\SetFigFont{6}{7.2}{rm}#2}}}
\put(3541,-3436){\makebox(0,0)[lb]{\smash{\SetFigFont{6}{7.2}{rm}#3}}}
\put(4741,-3436){\makebox(0,0)[lb]{\smash{\SetFigFont{6}{7.2}{rm}#4}}}
\put(7141,-3436){\makebox(0,0)[lb]{\smash{\SetFigFont{6}{7.2}{rm}#5}}}
\put(8341,-3436){\makebox(0,0)[lb]{\smash{\SetFigFont{6}{7.2}{rm}#6}}}

\end{picture}
}

\begingroup\makeatletter\ifx\SetFigFont\undefined
\def\x#1#2#3#4#5#6#7\relax{\def\x{#1#2#3#4#5#6}}%
\expandafter\x\fmtname xxxxxx\relax \def\y{splain}%
\ifx\x\y   
\gdef\SetFigFont#1#2#3{%
  \ifnum #1<17\tiny\else \ifnum #1<20\small\else
  \ifnum #1<24\normalsize\else \ifnum #1<29\large\else
  \ifnum #1<34\Large\else \ifnum #1<41\LARGE\else
     \huge\fi\fi\fi\fi\fi\fi
  \csname #3\endcsname}%
\else
\gdef\SetFigFont#1#2#3{\begingroup
  \count@#1\relax \ifnum 25<\count@\count@25\fi
  \def\x{\endgroup\@setsize\SetFigFont{#2pt}}%
  \expandafter\x
    \csname \romannumeral\the\count@ pt\expandafter\endcsname
    \csname @\romannumeral\the\count@ pt\endcsname
  \csname #3\endcsname}%
\fi
\fi\endgroup

\begingroup\makeatletter\ifx\SetFigFont\undefined
\def\x#1#2#3#4#5#6#7\relax{\def\x{#1#2#3#4#5#6}}%
\expandafter\x\fmtname xxxxxx\relax \def\y{splain}%
\ifx\x\y   
\gdef\SetFigFont#1#2#3{%
  \ifnum #1<17\tiny\else \ifnum #1<20\small\else
  \ifnum #1<24\normalsize\else \ifnum #1<29\large\else
  \ifnum #1<34\Large\else \ifnum #1<41\LARGE\else
     \huge\fi\fi\fi\fi\fi\fi
  \csname #3\endcsname}%
\else
\gdef\SetFigFont#1#2#3{\begingroup
  \count@#1\relax \ifnum 25<\count@\count@25\fi
  \def\x{\endgroup\@setsize\SetFigFont{#2pt}}%
  \expandafter\x
    \csname \romannumeral\the\count@ pt\expandafter\endcsname
    \csname @\romannumeral\the\count@ pt\endcsname
  \csname #3\endcsname}%
\fi
\fi\endgroup

\begingroup\makeatletter\ifx\SetFigFont\undefined
\def\x#1#2#3#4#5#6#7\relax{\def\x{#1#2#3#4#5#6}}%
\expandafter\x\fmtname xxxxxx\relax \def\y{splain}%
\ifx\x\y   
\gdef\SetFigFont#1#2#3{%
  \ifnum #1<17\tiny\else \ifnum #1<20\small\else
  \ifnum #1<24\normalsize\else \ifnum #1<29\large\else
  \ifnum #1<34\Large\else \ifnum #1<41\LARGE\else
     \huge\fi\fi\fi\fi\fi\fi
  \csname #3\endcsname}%
\else
\gdef\SetFigFont#1#2#3{\begingroup
  \count@#1\relax \ifnum 25<\count@\count@25\fi
  \def\x{\endgroup\@setsize\SetFigFont{#2pt}}%
  \expandafter\x
    \csname \romannumeral\the\count@ pt\expandafter\endcsname
    \csname @\romannumeral\the\count@ pt\endcsname
  \csname #3\endcsname}%
\fi
\fi\endgroup

\begingroup\makeatletter\ifx\SetFigFont\undefined
\def\x#1#2#3#4#5#6#7\relax{\def\x{#1#2#3#4#5#6}}%
\expandafter\x\fmtname xxxxxx\relax \def\y{splain}%
\ifx\x\y   
\gdef\SetFigFont#1#2#3{%
  \ifnum #1<17\tiny\else \ifnum #1<20\small\else
  \ifnum #1<24\normalsize\else \ifnum #1<29\large\else
  \ifnum #1<34\Large\else \ifnum #1<41\LARGE\else
     \huge\fi\fi\fi\fi\fi\fi
  \csname #3\endcsname}%
\else
\gdef\SetFigFont#1#2#3{\begingroup
  \count@#1\relax \ifnum 25<\count@\count@25\fi
  \def\x{\endgroup\@setsize\SetFigFont{#2pt}}%
  \expandafter\x
    \csname \romannumeral\the\count@ pt\expandafter\endcsname
    \csname @\romannumeral\the\count@ pt\endcsname
  \csname #3\endcsname}%
\fi
\fi\endgroup

\begingroup\makeatletter\ifx\SetFigFont\undefined
\def\x#1#2#3#4#5#6#7\relax{\def\x{#1#2#3#4#5#6}}%
\expandafter\x\fmtname xxxxxx\relax \def\y{splain}%
\ifx\x\y   
\gdef\SetFigFont#1#2#3{%
  \ifnum #1<17\tiny\else \ifnum #1<20\small\else
  \ifnum #1<24\normalsize\else \ifnum #1<29\large\else
  \ifnum #1<34\Large\else \ifnum #1<41\LARGE\else
     \huge\fi\fi\fi\fi\fi\fi
  \csname #3\endcsname}%
\else
\gdef\SetFigFont#1#2#3{\begingroup
  \count@#1\relax \ifnum 25<\count@\count@25\fi
  \def\x{\endgroup\@setsize\SetFigFont{#2pt}}%
  \expandafter\x
    \csname \romannumeral\the\count@ pt\expandafter\endcsname
    \csname @\romannumeral\the\count@ pt\endcsname
  \csname #3\endcsname}%
\fi
\fi\endgroup

\begingroup\makeatletter\ifx\SetFigFont\undefined
\def\x#1#2#3#4#5#6#7\relax{\def\x{#1#2#3#4#5#6}}%
\expandafter\x\fmtname xxxxxx\relax \def\y{splain}%
\ifx\x\y   
\gdef\SetFigFont#1#2#3{%
  \ifnum #1<17\tiny\else \ifnum #1<20\small\else
  \ifnum #1<24\normalsize\else \ifnum #1<29\large\else
  \ifnum #1<34\Large\else \ifnum #1<41\LARGE\else
     \huge\fi\fi\fi\fi\fi\fi
  \csname #3\endcsname}%
\else
\gdef\SetFigFont#1#2#3{\begingroup
  \count@#1\relax \ifnum 25<\count@\count@25\fi
  \def\x{\endgroup\@setsize\SetFigFont{#2pt}}%
  \expandafter\x
    \csname \romannumeral\the\count@ pt\expandafter\endcsname
    \csname @\romannumeral\the\count@ pt\endcsname
  \csname #3\endcsname}%
\fi
\fi\endgroup

\newcommand{\Div}[4]{\parbox[c]{1.8cm}{$
\begin{picture}(2751,2200)(6925,-5050)
\thinlines
\put(7201,-3961){\circle{336}}
\put(8401,-3961){\circle{336}}
\put(9253,-3109){\circle{336}}
\put(9253,-4813){\circle{336}}
\put(7369,-3961){\line( 1, 0){870}}
\put(8519,-3843){\line( 1, 1){615}}
\put(8519,-4079){\line( 1, -1){615}}
\put(7141,-3451){\makebox(0,0)[lb]{\smash{\SetFigFont{6}{7.2}{rm}#1}}}
\put(8341,-3511){\makebox(0,0)[lb]{\smash{\SetFigFont{6}{7.2}{rm}#2}}}
\put(9676,-3211){\makebox(0,0)[lb]{\smash{\SetFigFont{6}{7.2}{rm}#3}}}
\put(9676,-4936){\makebox(0,0)[lb]{\smash{\SetFigFont{6}{7.2}{rm}#4}}}
\end{picture}$}
}

\begingroup\makeatletter\ifx\SetFigFont\undefined
\def\x#1#2#3#4#5#6#7\relax{\def\x{#1#2#3#4#5#6}}%
\expandafter\x\fmtname xxxxxx\relax \def\y{splain}%
\ifx\x\y   
\gdef\SetFigFont#1#2#3{%
  \ifnum #1<17\tiny\else \ifnum #1<20\small\else
  \ifnum #1<24\normalsize\else \ifnum #1<29\large\else
  \ifnum #1<34\Large\else \ifnum #1<41\LARGE\else
     \huge\fi\fi\fi\fi\fi\fi
  \csname #3\endcsname}%
\else
\gdef\SetFigFont#1#2#3{\begingroup
  \count@#1\relax \ifnum 25<\count@\count@25\fi
  \def\x{\endgroup\@setsize\SetFigFont{#2pt}}%
  \expandafter\x
    \csname \romannumeral\the\count@ pt\expandafter\endcsname
    \csname @\romannumeral\the\count@ pt\endcsname
  \csname #3\endcsname}%
\fi
\fi\endgroup

\newcommand{\Dn}[5]{\parbox[c]{3.4cm}{$
\begin{picture}(5151,2200)(4525,-5050)
\thinlines
\put(4801,-3961){\circle{336}}
\put(6001,-3961){\circle{336}}
\put(8401,-3961){\circle{336}}
\put(9253,-3109){\circle{336}}
\put(9253,-4813){\circle{336}}
\put(7036,-3946){\circle*{30}}
\put(7186,-3946){\circle*{30}}
\put(7336,-3946){\circle*{30}}
\put(4969,-3961){\line( 1, 0){864}}
\put(6169,-3961){\line( 1, 0){600}}
\put(7603,-3961){\line( 1, 0){600}}
\put(8519,-3843){\line( 1, 1){615}}
\put(8519,-4079){\line( 1, -1){615}}
\put(4741,-3436){\makebox(0,0)[lb]{\smash{\SetFigFont{6}{7.2}{rm}\mbox{$#1$}}}}
\put(5941,-3436){\makebox(0,0)[lb]{\smash{\SetFigFont{6}{7.2}{rm}\mbox{$#2$}}}}
\put(8341,-3436){\makebox(0,0)[lb]{\smash{\SetFigFont{6}{7.2}{rm}\mbox{$#3$}}}}
\put(9676,-3211){\makebox(0,0)[lb]{\smash{\SetFigFont{6}{7.2}{rm}\mbox{$#4$}}}}
\put(9676,-4936){\makebox(0,0)[lb]{\smash{\SetFigFont{6}{7.2}{rm}\mbox{$#5$}}}}
\end{picture}$}
}
\begingroup\makeatletter\ifx\SetFigFont\undefined
\def\x#1#2#3#4#5#6#7\relax{\def\x{#1#2#3#4#5#6}}%
\expandafter\x\fmtname xxxxxx\relax \def\y{splain}%
\ifx\x\y   
\gdef\SetFigFont#1#2#3{%
  \ifnum #1<17\tiny\else \ifnum #1<20\small\else
  \ifnum #1<24\normalsize\else \ifnum #1<29\large\else
  \ifnum #1<34\Large\else \ifnum #1<41\LARGE\else
     \huge\fi\fi\fi\fi\fi\fi
  \csname #3\endcsname}%
\else
\gdef\SetFigFont#1#2#3{\begingroup
  \count@#1\relax \ifnum 25<\count@\count@25\fi
  \def\x{\endgroup\@setsize\SetFigFont{#2pt}}%
  \expandafter\x
    \csname \romannumeral\the\count@ pt\expandafter\endcsname
    \csname @\romannumeral\the\count@ pt\endcsname
  \csname #3\endcsname}%
\fi
\fi\endgroup

\begingroup\makeatletter\ifx\SetFigFont\undefined
\def\x#1#2#3#4#5#6#7\relax{\def\x{#1#2#3#4#5#6}}%
\expandafter\x\fmtname xxxxxx\relax \def\y{splain}%
\ifx\x\y   
\gdef\SetFigFont#1#2#3{%
  \ifnum #1<17\tiny\else \ifnum #1<20\small\else
  \ifnum #1<24\normalsize\else \ifnum #1<29\large\else
  \ifnum #1<34\Large\else \ifnum #1<41\LARGE\else
     \huge\fi\fi\fi\fi\fi\fi
  \csname #3\endcsname}%
\else
\gdef\SetFigFont#1#2#3{\begingroup
  \count@#1\relax \ifnum 25<\count@\count@25\fi
  \def\x{\endgroup\@setsize\SetFigFont{#2pt}}%
  \expandafter\x
    \csname \romannumeral\the\count@ pt\expandafter\endcsname
    \csname @\romannumeral\the\count@ pt\endcsname
  \csname #3\endcsname}%
\fi
\fi\endgroup

\begingroup\makeatletter\ifx\SetFigFont\undefined
\def\x#1#2#3#4#5#6#7\relax{\def\x{#1#2#3#4#5#6}}%
\expandafter\x\fmtname xxxxxx\relax \def\y{splain}%
\ifx\x\y   
\gdef\SetFigFont#1#2#3{%
  \ifnum #1<17\tiny\else \ifnum #1<20\small\else
  \ifnum #1<24\normalsize\else \ifnum #1<29\large\else
  \ifnum #1<34\Large\else \ifnum #1<41\LARGE\else
     \huge\fi\fi\fi\fi\fi\fi
  \csname #3\endcsname}%
\else
\gdef\SetFigFont#1#2#3{\begingroup
  \count@#1\relax \ifnum 25<\count@\count@25\fi
  \def\x{\endgroup\@setsize\SetFigFont{#2pt}}%
  \expandafter\x
    \csname \romannumeral\the\count@ pt\expandafter\endcsname
    \csname @\romannumeral\the\count@ pt\endcsname
  \csname #3\endcsname}%
\fi
\fi\endgroup

\newcommand{\Dnii}[7]{\parbox[c]{4.8cm}{$
\begin{picture}(7551,2200)(4525,-5050)
\thinlines
\put(4801,-3961){\circle{336}}
\put(6001,-3961){\circle{336}}
\put(7201,-3961){\circle{336}}
\put(8401,-3961){\circle{336}}
\put(10801,-3961){\circle{336}}
\put(11653,-3109){\circle{336}}
\put(11653,-4813){\circle{336}}
\put(9436,-3946){\circle*{30}}
\put(9586,-3946){\circle*{30}}
\put(9736,-3946){\circle*{30}}
\put(4969,-3961){\line( 1, 0){864}}
\put(6169,-3961){\line( 1, 0){864}}
\put(7369,-3961){\line( 1, 0){864}}
\put(8569,-3961){\line( 1, 0){600}}
\put(10003,-3961){\line( 1, 0){600}}
\put(10919,-3843){\line( 1, 1){615}}
\put(10919,-4079){\line( 1, -1){615}}
\put(4741,-3436){\makebox(0,0)[lb]{\smash{\SetFigFont{6}{7.2}{rm}#1}}}
\put(5941,-3436){\makebox(0,0)[lb]{\smash{\SetFigFont{6}{7.2}{rm}#2}}}
\put(7141,-3436){\makebox(0,0)[lb]{\smash{\SetFigFont{6}{7.2}{rm}#3}}}
\put(8341,-3436){\makebox(0,0)[lb]{\smash{\SetFigFont{6}{7.2}{rm}#4}}}
\put(10741,-3436){\makebox(0,0)[lb]{\smash{\SetFigFont{6}{7.2}{rm}#5}}}
\put(12076,-3211){\makebox(0,0)[lb]{\smash{\SetFigFont{6}{7.2}{rm}#6}}}
\put(12076,-4936){\makebox(0,0)[lb]{\smash{\SetFigFont{6}{7.2}{rm}#7}}}
\end{picture}$}
}
\begingroup\makeatletter\ifx\SetFigFont\undefined
\def\x#1#2#3#4#5#6#7\relax{\def\x{#1#2#3#4#5#6}}%
\expandafter\x\fmtname xxxxxx\relax \def\y{splain}%
\ifx\x\y   
\gdef\SetFigFont#1#2#3{%
  \ifnum #1<17\tiny\else \ifnum #1<20\small\else
  \ifnum #1<24\normalsize\else \ifnum #1<29\large\else
  \ifnum #1<34\Large\else \ifnum #1<41\LARGE\else
     \huge\fi\fi\fi\fi\fi\fi
  \csname #3\endcsname}%
\else
\gdef\SetFigFont#1#2#3{\begingroup
  \count@#1\relax \ifnum 25<\count@\count@25\fi
  \def\x{\endgroup\@setsize\SetFigFont{#2pt}}%
  \expandafter\x
    \csname \romannumeral\the\count@ pt\expandafter\endcsname
    \csname @\romannumeral\the\count@ pt\endcsname
  \csname #3\endcsname}%
\fi
\fi\endgroup

\begingroup\makeatletter\ifx\SetFigFont\undefined
\def\x#1#2#3#4#5#6#7\relax{\def\x{#1#2#3#4#5#6}}%
\expandafter\x\fmtname xxxxxx\relax \def\y{splain}%
\ifx\x\y   
\gdef\SetFigFont#1#2#3{%
  \ifnum #1<17\tiny\else \ifnum #1<20\small\else
  \ifnum #1<24\normalsize\else \ifnum #1<29\large\else
  \ifnum #1<34\Large\else \ifnum #1<41\LARGE\else
     \huge\fi\fi\fi\fi\fi\fi
  \csname #3\endcsname}%
\else
\gdef\SetFigFont#1#2#3{\begingroup
  \count@#1\relax \ifnum 25<\count@\count@25\fi
  \def\x{\endgroup\@setsize\SetFigFont{#2pt}}%
  \expandafter\x
    \csname \romannumeral\the\count@ pt\expandafter\endcsname
    \csname @\romannumeral\the\count@ pt\endcsname
  \csname #3\endcsname}%
\fi
\fi\endgroup

\begingroup\makeatletter\ifx\SetFigFont\undefined
\def\x#1#2#3#4#5#6#7\relax{\def\x{#1#2#3#4#5#6}}%
\expandafter\x\fmtname xxxxxx\relax \def\y{splain}%
\ifx\x\y   
\gdef\SetFigFont#1#2#3{%
  \ifnum #1<17\tiny\else \ifnum #1<20\small\else
  \ifnum #1<24\normalsize\else \ifnum #1<29\large\else
  \ifnum #1<34\Large\else \ifnum #1<41\LARGE\else
     \huge\fi\fi\fi\fi\fi\fi
  \csname #3\endcsname}%
\else
\gdef\SetFigFont#1#2#3{\begingroup
  \count@#1\relax \ifnum 25<\count@\count@25\fi
  \def\x{\endgroup\@setsize\SetFigFont{#2pt}}%
  \expandafter\x
    \csname \romannumeral\the\count@ pt\expandafter\endcsname
    \csname @\romannumeral\the\count@ pt\endcsname
  \csname #3\endcsname}%
\fi
\fi\endgroup

\newcommand{\Fiv}[4]{
\begin{picture}(3952,895)(1025,-4136)
\thinlines
\put(1201,-3961){\circle{336}}
\put(2401,-3961){\circle{336}}
\put(3601,-3961){\circle*{336}}
\put(4801,-3961){\circle*{336}}
\put(1369,-3961){\line( 1, 0){864}}
\put(2569,-3911){\line( 1, 0){864}}
\put(2569,-4011){\line( 1, 0){864}}
\put(3769,-3961){\line( 1, 0){864}}
\put(1141,-3436){\makebox(0,0)[lb]{\smash{\SetFigFont{6}{7.2}{rm}\mbox{$#1$}}}}
\put(2341,-3436){\makebox(0,0)[lb]{\smash{\SetFigFont{6}{7.2}{rm}\mbox{$#2$}}}}
\put(3541,-3436){\makebox(0,0)[lb]{\smash{\SetFigFont{6}{7.2}{rm}\mbox{$#3$}}}}
\put(4741,-3436){\makebox(0,0)[lb]{\smash{\SetFigFont{6}{7.2}{rm}\mbox{$#4$}}}}
\end{picture}
}
\begingroup\makeatletter\ifx\SetFigFont\undefined
\def\x#1#2#3#4#5#6#7\relax{\def\x{#1#2#3#4#5#6}}%
\expandafter\x\fmtname xxxxxx\relax \def\y{splain}%
\ifx\x\y   
\gdef\SetFigFont#1#2#3{%
  \ifnum #1<17\tiny\else \ifnum #1<20\small\else
  \ifnum #1<24\normalsize\else \ifnum #1<29\large\else
  \ifnum #1<34\Large\else \ifnum #1<41\LARGE\else
     \huge\fi\fi\fi\fi\fi\fi
  \csname #3\endcsname}%
\else
\gdef\SetFigFont#1#2#3{\begingroup
  \count@#1\relax \ifnum 25<\count@\count@25\fi
  \def\x{\endgroup\@setsize\SetFigFont{#2pt}}%
  \expandafter\x
    \csname \romannumeral\the\count@ pt\expandafter\endcsname
    \csname @\romannumeral\the\count@ pt\endcsname
  \csname #3\endcsname}%
\fi
\fi\endgroup

\newcommand{\Gii}[2]{
\begin{picture}(1552,895)(1025,-4136)
\thinlines
\put(1201,-3961){\circle*{336}}
\put(2401,-3961){\circle{336}}
\put(1339,-3861){\line( 1, 0){924}}
\put(1369,-3961){\line( 1, 0){864}}
\put(1339,-4061){\line( 1, 0){924}}
\put(1141,-3436){\makebox(0,0)[lb]{\smash{\SetFigFont{6}{7.2}{rm}\mbox{$#1$}}}}
\put(2341,-3436){\makebox(0,0)[lb]{\smash{\SetFigFont{6}{7.2}{rm}\mbox{$#2$}}}}
\end{picture}
}
\begingroup\makeatletter\ifx\SetFigFont\undefined
\def\x#1#2#3#4#5#6#7\relax{\def\x{#1#2#3#4#5#6}}%
\expandafter\x\fmtname xxxxxx\relax \def\y{splain}%
\ifx\x\y   
\gdef\SetFigFont#1#2#3{%
  \ifnum #1<17\tiny\else \ifnum #1<20\small\else
  \ifnum #1<24\normalsize\else \ifnum #1<29\large\else
  \ifnum #1<34\Large\else \ifnum #1<41\LARGE\else
     \huge\fi\fi\fi\fi\fi\fi
  \csname #3\endcsname}%
\else
\gdef\SetFigFont#1#2#3{\begingroup
  \count@#1\relax \ifnum 25<\count@\count@25\fi
  \def\x{\endgroup\@setsize\SetFigFont{#2pt}}%
  \expandafter\x
    \csname \romannumeral\the\count@ pt\expandafter\endcsname
    \csname @\romannumeral\the\count@ pt\endcsname
  \csname #3\endcsname}%
\fi
\fi\endgroup

\newcommand{\Evi}[6]{\raisebox{2.4mm}{\parbox[c]{3.4cm}{$
\begin{picture}(5300,2200)(950,-5400)
\thinlines
\put(1201,-3961){\circle{336}}
\put(2401,-3961){\circle{336}}
\put(3601,-3961){\circle{336}}
\put(4801,-3961){\circle{336}}
\put(6001,-3961){\circle{336}}
\put(3601,-5161){\circle{336}}
\put(1369,-3961){\line( 1, 0){864}}
\put(2569,-3961){\line( 1, 0){864}}
\put(3769,-3961){\line( 1, 0){864}}
\put(4969,-3961){\line( 1, 0){864}}
\put(3601,-4129){\line( 0, -1){864}}
\put(1141,-3436){\makebox(0,0)[lb]{\smash{\SetFigFont{6}{7.2}{rm}\mbox{$#1$}}}}
\put(2341,-3436){\makebox(0,0)[lb]{\smash{\SetFigFont{6}{7.2}{rm}\mbox{$#3$}}}}
\put(3541,-3436){\makebox(0,0)[lb]{\smash{\SetFigFont{6}{7.2}{rm}\mbox{$#4$}}}}
\put(4741,-3436){\makebox(0,0)[lb]{\smash{\SetFigFont{6}{7.2}{rm}\mbox{$#5$}}}}
\put(5941,-3436){\makebox(0,0)[lb]{\smash{\SetFigFont{6}{7.2}{rm}\mbox{$#6$}}}}
\put(4024,-5263){\makebox(0,0)[lb]{\smash{\SetFigFont{6}{7.2}{rm}\mbox{$#2$}}}}
\end{picture}$}}
}

\begingroup\makeatletter\ifx\SetFigFont\undefined
\def\x#1#2#3#4#5#6#7\relax{\def\x{#1#2#3#4#5#6}}%
\expandafter\x\fmtname xxxxxx\relax \def\y{splain}%
\ifx\x\y   
\gdef\SetFigFont#1#2#3{%
  \ifnum #1<17\tiny\else \ifnum #1<20\small\else
  \ifnum #1<24\normalsize\else \ifnum #1<29\large\else
  \ifnum #1<34\Large\else \ifnum #1<41\LARGE\else
     \huge\fi\fi\fi\fi\fi\fi
  \csname #3\endcsname}%
\else
\gdef\SetFigFont#1#2#3{\begingroup
  \count@#1\relax \ifnum 25<\count@\count@25\fi
  \def\x{\endgroup\@setsize\SetFigFont{#2pt}}%
  \expandafter\x
    \csname \romannumeral\the\count@ pt\expandafter\endcsname
    \csname @\romannumeral\the\count@ pt\endcsname
  \csname #3\endcsname}%
\fi
\fi\endgroup

\begingroup\makeatletter\ifx\SetFigFont\undefined
\def\x#1#2#3#4#5#6#7\relax{\def\x{#1#2#3#4#5#6}}%
\expandafter\x\fmtname xxxxxx\relax \def\y{splain}%
\ifx\x\y   
\gdef\SetFigFont#1#2#3{%
  \ifnum #1<17\tiny\else \ifnum #1<20\small\else
  \ifnum #1<24\normalsize\else \ifnum #1<29\large\else
  \ifnum #1<34\Large\else \ifnum #1<41\LARGE\else
     \huge\fi\fi\fi\fi\fi\fi
  \csname #3\endcsname}%
\else
\gdef\SetFigFont#1#2#3{\begingroup
  \count@#1\relax \ifnum 25<\count@\count@25\fi
  \def\x{\endgroup\@setsize\SetFigFont{#2pt}}%
  \expandafter\x
    \csname \romannumeral\the\count@ pt\expandafter\endcsname
    \csname @\romannumeral\the\count@ pt\endcsname
  \csname #3\endcsname}%
\fi
\fi\endgroup

\title{Representations of compact simple Lie groups\\
whose orbit space has non-empty boundary}
\author{Claudio Gorodski}
\thanks{The first author acknowledges partial financial
  support from  CNPq (grant 302882/2017-0) and FAPESP (grant 16/23746-6).}

\address{Instituto de Matem\'atica e Estat\'\i stica, Universidade de 
  S\~ao Paulo, Rua do Mat\~ao, 1010, S\~ao Paulo, SP 05508-090, Brazil}
\email{gorodski@ime.usp.br}

\author{Andreas Kollross}
\address{Universität Stuttgart, Institut für Geometrie und Topologie,
Pfaffenwaldring~57, 70569~Stuttgart, Germany}

\email{kollross@mathematik.uni-stuttgart.de}

\author{Burkhard Wilking}

\address{Mathematisches Institut, M\"unster University, Einsteinstrasse~62,
  48149~M\"unster, Germany}
\email{wilking@uni-muenster.de}

\date{\today}


\begin{document}

\maketitle

\begin{abstract}
  We provide detailed calculations for the classification of representations
  of compact simple Lie groups with non-empty boundary in the 
orbit space, first announced in the paper~\cite{GKW} by the same authors. 
\end{abstract}

\section{Introduction}

In~\cite{GKW} a proof of the following theorem is sketched. The purpose 
of the present text is to provide the detailed calculations. We refer
the reader to~\cite{GKW} for the relevant context and 
necessary definitions.  

\begin{thm}\label{classif}
  The representations $V$ of
compact connected simple Lie groups $G$ with non-empty boundary
in the orbit space are listed in Tables~1 and~2,
up to a trivial component and up to an outer automorphism.
In the irreducible case (Table~1),
we also indicate the kernel of the representation in those cases
in which it is non-trivial, the effective
principal isotropy group,
and whether the representation is polar, toric
or quaternion-toric.

\begin{table}[t]
\rm
\[ \begin{array}{|c|c|c|c|c|}
  \hline
  G & \textsl{Kernel} & V & \textsl{Property} & \textsl{Effective p.i.g.} \\
  \hline
  \SU2 &\mbox{---}& \C^2 & \textrm{polar} & 1 \\
  \hline
   \multirow{2}{*}{\SO3}
  &\multirow{2}{*}{\mbox{---}}& \R^3 &\multirow{2}{*}{polar}&\sf T^1\\
  &&\mathrm{S}^2_0\R^3=\R^5 && \sf \Z_2^2\\
  \hline
  \multirow{3}{*}{\shortstack{\SU n\\ ($n\geq3$)}} &\mbox{---}&\C^n & \multirow{2}{*}{\textrm{polar}}&\SU{n-1}\\
 &\Z_n&\mathrm{Ad}&&\sf T^{n-1} \\ \cline{4-4}
  &\textrm{$\{\pm1\}$ if $n$ is even} &\mathrm{S}^2\C^n& \textrm{toric}&\Z_2^{n-1}/\mathrm{ker}\\
  \hline
    \multirow{2}{*}{\shortstack{\SU n\\ ($n\geq5$)}}
   &\multirow{2}{*}{\textrm{$\{\pm1\}$ if $n$ is even}}
   &\multirow{2}{*}{$\Lambda^2\C^n$}
   &\parbox[t]{.8in}{polar ($n$~odd)}
   &\multirow{2}{*}{$\SU2^{\lfloor\frac n2\rfloor}/\mathrm{ker}$}\\
   &&&\parbox[t]{.8in}{toric ($n$~even)}&\\
  \hline
  \SU6 &\mbox{---}& \Lambda^3\C^6=\Q^{10} & \textrm{q-toric}&\sf T^2\\ \hline
  \SU8& \Z_4 & [\Lambda^4\C^8]_{\mathbb R} & \textrm{polar} & \Z_2^7 \\ \hline
  \multirow{3}{*}{\shortstack{\SO n\\ ($n\geq5$)}}  &\mbox{---}& \R^n &
  \multirow{3}{*}{polar}&\Spin{n-1}\\ \cline{2-2}
  &\multirow{2}{*}{\textrm{$\{\pm1\}$ if $n$ is even}}&\Lambda^2\R^n=\mathrm{Ad}&&\sf T^{\lfloor\frac n2\rfloor} \\
  &&\mathrm{S}^2_0\R^n&&\Z_2^{n-1}/\mathrm{ker}\\  \hline
  \Spin7&\mbox{---}&\R^8\;\textrm{(spin)}& \textrm{polar} &\G \\ \hline
  \Spin8 &\Z_2& \R^8_\pm\; \textrm{(half-spin)} &\textrm{polar} & \mathsf{Spin}'(7) \\  \hline
  \Spin9 & \mbox{---} &\R^{16}\;\textrm{(spin)} & \textrm{polar} & \Spin7 \\ \hline
  \Spin{10} & \mbox{---} &\C^{16}_\pm\; \textrm{(half-spin)} & \textrm{polar} & \SU4 \\ \hline
  \Spin{11} & \mbox{---}& \Q^{16}\;\textrm{(spin)} & \mbox{---} & 1 \\ \hline
  \Spin{12} & \Z_2 & \Q^{16}_\pm\;\textrm{(half-spin)} & \textrm{q-toric} & \SP1^3\\ \hline
    \Spin{16} & \Z_2 & \R^{128}_\pm\;\textrm{(half-spin)} & \textrm{polar} & \Z_2^8\\ \hline
  \multirow{3}{*}{\shortstack{\SP n\\ ($n\geq3$)}} &\mbox{---}& \C^{2n}=\Q^n & \multirow{3}{*}{polar}&\SP{n-1} \\  \cline{2-2}
  &\multirow{2}{*}{$\{\pm1\}$}&[\mathrm{S}^2\C^{2n}]_{\mathbb R}=\mathrm{Ad}&&\sf T^n \\
  &&[\Lambda^2_0\C^{2n}]_{\mathbb R} &&\SP1^n/\{\pm1\} \\
  \hline
  \SP3 &\mbox{---}& \Lambda^3_0\C^6=\Q^7 &\textrm{q-toric}&\Z_2^2 \\ \hline
  \SP4 & \{\pm1\} & [\Lambda^4_0\C^8]_{\mathbb R} & \textrm{polar}&\Z_2^6 \\ \hline
  \multirow{2}{*}{\G} &\multirow{2}{*}{\mbox{---}}& \R^7 & \multirow{2}{*}{polar} &\SU3 \\
  &&\mathrm{Ad} &&\sf T^2\\ \hline
  \multirow{2}{*}{\F} &\multirow{2}{*}{\mbox{---}}& \R^{26}
  &\multirow{2}{*}{polar}&\Spin8\\
  &&\mathrm{Ad} &&\sf T^4\\ \hline
  \E6 &\mbox{---}& \C^{27} & \textrm{toric} & \Spin8\\ \hline
  \E6& \Z_3 &\mathrm{Ad} & \textrm{polar} & \sf T^6\\ \hline
  \E7 &\mbox{---}& \Q^{28} & \textrm{q-toric} & \Spin8\\ \hline
    \E7 & \Z_2 &\mathrm{Ad} &  \textrm{polar} & \sf T^7 \\\hline
  \E8 &\mbox{---}& \mathrm{Ad} & \textrm{polar} & \sf T^8\\
   \hline
\end{array} \]
\smallskip
\begin{center}
  \sc Table~1: Irreducible representations of compact simple
  Lie groups with non-empty boundary in the orbit space.
\end{center}
\end{table}

\begin{table}[t]
\[ \begin{array}{|c|c|c|}
  \hline
  \multirow{2}{*}{\SU n} & k\,\C^n & 2\leq k\leq n-1 \\
  & \C^n\oplus \Lambda^2\C^n  & n\geq5 \\
  \hline
  \multirow{2}{*}{\SU4} & k\,\R^6\oplus\ell\,\C^4 &  2\leq k+\ell\leq 3\\
  & \R^6\oplus\mathrm{Ad} & - \\
  \hline
  \multirow{2}{*}{\SO n} & k\,\R^n & 2\leq k\leq n-1\\
  & \R^n\oplus\mathrm{Ad} & n\geq5 \\
  \hline
  \SP2 & \Q^2\oplus\R^5 & -\\
  \hline
  \Spin7 & k\,\R^7\oplus\ell\,\R^8 & 2\leq k+\ell\leq 4 \\
  \hline
  \Spin8 & k\,\R^8\oplus\ell\,\R^8_+\oplus m\,\R^8_- & 2\leq k+\ell+m\leq5\\
  \hline
  \multirow{3}{*}{\Spin9} & k\, \R^{16} & 2\leq k\leq3\\
  & \R^{16}\oplus k\,\R^9 & 1\leq k\leq4\\
  &2\R^{16}\oplus k\,\R^9 & 0\leq k\leq2\\
  \hline
  \Spin{10} & \C^{16}\oplus k\,\R^{10} & 1\leq k\leq 3\\
  \hline
  \Spin{12} & \Q^{16}\oplus\R^{12} & -\\
  \hline
  \multirow{2}{*}{\SP n} & k\,\C^{2n} & 2\leq k\leq n\\
  & \C^{2n} \oplus[\Lambda^2_0\C^{2n}]_{\mathbb R} & n\geq3 \\ \hline
  \SP3 & 2\,[\Lambda^2_0\C^6]_{\mathbb R} & - \\ \hline
  \G & k\,\R^7 & 2\leq k\leq3\\
  \hline
  \F & 2\,\R^{26} & -\\
  \hline
\end{array} \]
\smallskip
\begin{center}
  \sc Table~2: Reducible representations of compact simple
  Lie groups with non-empty boundary in the orbit space.
\end{center}
\footnotesize{\rm In case of $\Spin8$,
  the prime in $\mathsf{Spin}'(7)$ refers to a nonstandard $\Spin7$-subgroup;
  in case of $\Spin n$, $\mathrm{S}^2_0\R^n=\mathrm{S}^2\R^n\ominus\R$;
  in case of $\SP n$,
  $\Lambda^k_0\C^{2n}=\Lambda^k\C^{2n}\ominus\Lambda^{k-2}\C^{2n}$;
  and $[V]_{\mathbb R}$ denotes a real form of~$V$.}
\end{table}
\end{thm}

\section{The irreducible case}

Assume $G$ is a simple Lie group, 
$\rho:G\to \OG V$ is irreducible and $X=V/G$ has non-empty 
boundary. The first step is checking representations up to 
dimension~$\mathcal L_G$.  The real irreducible representations 
of compact Lie groups are either real forms or realifications 
of complex irreducible representations; which case does happen can be read off
the highest weight vector of the complex irreducible
representation~\cite[App.~A]{GT}.  
Also, Weyl's dimension formula gives the 
complex dimension of an irreducible representation in terms of the 
highest weight vector, and for each simple Lie group there are 
certain obvious 
monotonous dependences of the dimension in terms of coefficients of the highest
weight vector~\cite{On}.  

Polar representations of compact Lie groups have non-empty 
boundary in the orbit space, as well as toric and q-toric
ones~\cite[\S~1.3]{GKW}. We shall refer to these as the
standard examples of representations with non-empty boundary
in the orbit space. 
The main task facing us here is to look for non-standard examples
for simple Lie groups.
As we shall see, it turns out there is only one:
the spin representation of $\Spin{11}$. 

\subsection{The case of $\SU n$}

The standard examples are listed on Table~1. We shall refer
to $\C^n$, $\mathrm{Ad}$, $S^2\C^n$ and $\Lambda^2\C^n$ as the
basic examples, and to $(\SO3,\R^5)$, $(\SO6,S_0^2\R^6)$,
$(\SO6,\R^6)$,
$(\SU6,\Lambda^3\C^6)$, $(\SU8,[\Lambda^4\C^8]_{\mathbb R})$
as the special examples, as these occur only for specific values of~$n$. 

The case of $\SU2$ is studied in~\cite[\S10, \S11]{GL}, where it is shown
that $\rho$ is polar, so $V=\C^2$, $\R^3$ or $\R^5$. 

In the case of $\SU3$, we list the representations of dimension 
up to~$\mathcal L_G=24$ and get only $S^3\C^3$ of dimension~$20$,
in addition to the standard representations. 
We are going to check that there are no
nice involutions for $S^3\C^3$,
which is a contradiction to~\cite[Lemma~3.1]{GKW}. 
It is easy to see that any involution in $\SU3/Z(\SU3)$ is conjugate
to one of $\omega\left(\begin{smallmatrix} 1&0&0\\0&-1&0\\0&0&-1\end{smallmatrix}\right)$, where $\omega^3=1$,  or $\omega\left(\begin{smallmatrix} -1&0&0\\0&1&0\\0&0&1\end{smallmatrix}\right)$, where $\omega^3=-1$. Both elements have the same fixed point set
in $S^3(\C^3)$, which has codimension~$12$, larger than 
$4+\dim\C P^2=8$, which is the bound allowed for a nice involution of 
$\SU3$. 

In the case of $\SU4$, 
$\mathcal L_G=40$ and,
besides the standard representations, 
we only get
\begin{equation}\label{a3}
  \Aiii11{}
\end{equation}
of dimension exactly~$40$.
Any involution in $\SU4/Z(\SU4)$ is conjugate
to one of $\omega\mathrm{diag}(-1_2,1_2)$, where $\omega^4=1$,
or $\omega\mathrm{diag}(-1,1_3)$, where $\omega^4=-1$. 
We can determine their fixed point set in~(\ref{a3})  
by looking at the decomposition
\[ \Aiii1{}{}{}\otimes\Aiii{}1{}=\Aiii11{}\oplus\Aiii{}{}1. \]
The left hand-side is $\C^4\otimes\Lambda^2\C^4$; we readily see that 
only the first kind of involutions has a nonzero fixed point subspace,
and that is of dimension $12$. Therefore their fixed point subspace 
in~(\ref{a3}) has dimension~$10$, thus codimension~$30$, 
which is bigger than $4+\dim\SU4/\mathsf{S}(\U2\times\U2)=12$, so that 
there are no nice involutions. 

In the case of $\SU5$,  $\mathcal L_G=60$ and we get no non-standard
representations within the dimension bound. 

In the case of $\SU6$,  $\mathcal L_G=84$ and we get no non-standard
representations within the dimension bound. 

In the case of $\SU7$,  $\mathcal L_G=112$ and the only non-standard
representation
within the dimension bound is $\Lambda^3\C^7$, of dimension $70$. 
We show that this 
representation admits no nice involution and hence has empty boundary. 
Every involution $\sigma\in\SU7/Z(\SU7)$ is conjugate to one of  
$\omega\mathrm{diag}(1,-1_6)$, $\omega\mathrm{diag}(-1_2,1_5)$, 
$\omega\mathrm{diag}(1_3,-1_4)$, where $\omega^7=1$,
or $\omega\mathrm{diag}(-1_k,1_\ell)$ where $k$ is odd and $\omega^7=-1$. 
The latter involution has no non-zero fixed points in $\Lambda^3\C^7$. 
The first three involutions have non-zero fixed points only if $\omega=1$, 
in which case they span a subspace of dimension  $32$, $30$ and $38$, resp.,
so the codimension is $38$, $40$ and $32$, which is bigger than the 
bound $4+\dim G/G^\sigma$ allowed for a nice involution, namely, $16$, 
$24$ and $28$. 

In the case of $\SU8$,  $\mathcal L_G=144$ and the only non-polar
representation
within the dimension bound is $\Lambda^3\C^8$,
of dimension $112$. 
We show that $\Lambda^3\C^8$
admits no nice involution and hence has empty boundary. 
Every involution with non-zero subspace of 
fixed points in that representation
is conjugate to one of 
$\omega\mathrm{diag}(-1_2,1_6)$, $\omega\mathrm{diag}(-1_4,1_4)$. 
$\omega\mathrm{diag}(-1_6,1_2)$, where $\omega^8=1$. 
Their dimensions are $2\binom 63+6=52$, $2(\binom 43+4\binom 42)=56$
$4\binom 62=60$, hence the codimensions are $60$, $56$, $52$, which is larger
than  $4+\dim G/G^\sigma=28$, $36$, $28$, so there cannot be nice
involutions. 

In the case of $\SU9$,  $\mathcal L_G=180$ and the only non-standard
representation within the dimension bound is $\Lambda^3\C^9$,
of dimension~$168$.
We repeat the argument above: the involutions that need to be considered are 
$\omega\mathrm{diag}(-1_k,1_{9-k})$, where $\omega^3=1$ and 
$k=2$, $4$, $6$, $8$, or
$\omega\mathrm{diag}(-1_k,1_{9-k})$, where $\omega^3=-1$ and 
$k=1$, $3$, $5$, $7$. The subspace of fixed points 
has codimension $84$, $88$, $76$, $112$, $112$, $76$, $88$, $84$, resp.,
which is bigger than $4+\dim G/G^\sigma=32$, $44$, $40$, $20$, $20$,
$40$, $44$, $32$, so there cannot be nice involutions. 

Consider $\Lambda^3\C^n$ for $n\geq10$. It is of complex type
and has dimension $\frac13n(n-1)(n-2)>2n^2+2n=\mathcal L_G$.

We refer to~\cite[Table~5]{OV} for the claims below about 
the dimensions of certain representations. 

Consider 
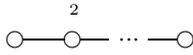
\[ \An{}2{}{} \]
It has dimension $\frac23\binom n2\binom{n+1}2>2n^2+2n$ for $n\geq5$. 

Consider 
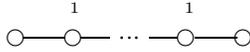
\[ \Ani{}11{} \]
It has dimension $\frac{n-3}{n-1}\binom n2\binom{n+1}2>2n^2+2n$ for 
$n\geq7$. 

Consider 
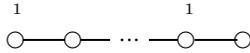
\[ \Ani1{}1{} \]
It has dimension $\frac{n(n-2)(n+1)^2}{n-1}>2n^2+2n$ for 
$n\geq4$.

Consider 
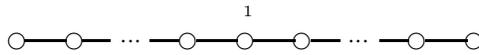
\[ \Anii{}{}{}1{}{}{} \] 
for $n=4k$. It is of real type and has dimension 
$\binom{4k}{2k}>2(4k)^2+2(4k)$ for $k\geq3$. 

This is enough to deduce that there will be no other 
irreducible representations within the dimension bound.

\subsection{The case of $\Spin{2n+1}$}

Here $\mathcal L_G=4n^2+10n+4$ for $n\geq1$.
The basic examples 
are the vector representation
$\R^{2n+1}$, the adjoint representation $\mathfrak{so}(2n+1)=\Lambda^2\R^{2n+1}$, 
the reduced symmetric square $S^2_0(\R^{2n+1})$, and the special examples
are the spin representations
of $\Spin5=\SP2$, $\Spin7$ and $\Spin9$ on $\Q^2$,
$\R^8$ and $\R^{16}$, respectively.

We first study the spin representation $\Delta_{2n+1}$ of $\Spin{2n+1}$.
\[ \mbox{$\Delta_{2n+1}$ is of}\quad\begin{array}{ll}\mbox{real type}&\mbox{if $n=4k$ or $4k+3$,}\\
\mbox{quaternionic type}&\mbox {if $n=4k+1$ or $4k+2$.}\end{array}. \]
Also, $\dim_{\mathbb C}\Delta_{2n+1}^{\mathbb C}=2^n$ and
\begin{equation}\label{d} 
\dim_{\mathbb R}\Delta_{2n+1}=\left\{\begin{array}{ll}2^n&\mbox{if $n=4k$ or $4k+3$,}\\
  2^{n+1}&\mbox {if $n=4k+1$ or $4k+2$.}\end{array}.\right. 
\end{equation}
The only non-polar cases
within the dimension bound are $\Delta_{11}$, $\Delta_{13}$, $\Delta_{15}$
and~$\Delta_{17}$. In~\cite[\S~8.1]{GKW}
it is shown that the orbit space of $\Delta_{11}$ has
non-empty boundary. In the following we apply Borel-de Siebenthal
theory~\cite[Theorem~8.10.8]{Wo}
to describe the involutions of the groups and prove that they do not
act as nice involutions on the other three
representations. It is relevant to point out that
Borel de-Siebenthal theory yields representatives for
the involutions modulo the center, but we still need to take into
account the action of the center under the representation
in order to compute the dimension of the fixed point subspace. 

The roots of $\Spin{2n+1}$ are $\pm\theta_i\pm\theta_j$ for $1\leq i<j\leq n$
and $\pm\theta_i$ for $1\leq i \leq n$. The system of simple roots
is $\Pi=\{\alpha_1,\ldots,\alpha_n\}$, where $\theta_i=\theta_i-\theta_{i+1}$
for $i=1,\ldots,n-1$ and $\alpha_n=\theta_n$. The highest root is
$\beta=\theta_1+\theta_2=\alpha_1+2(\alpha_2+\cdots+\alpha_n)$. 
Consider the vertices of the Cartan polyhedron
$v_0=0$, $v_1,\ldots,v_n\in\sqrt{-1}\Lt$, where $\Lt$ is the Lie algebra
of the maximal torus, given by
\[ \alpha_i(v_j)=\left\{\begin{array}{ll}1/m_i,&\mbox{if $i=j$}\\

0,&\mbox{if $i\neq j$}.\end{array}\right. \]
Here the $m_i$ are the coefficients of the highest root.
Then
\begin{eqnarray*}
  v_1&=&(1,0,\ldots,0),\\
  v_2&=&(\x,\x,0,\ldots,0),\\
  v_3&=&(\x,\x,\x,0,\ldots,0),\\
\vdots\\
  v_n&=&(\x,\ldots,\x),
  \end{eqnarray*}
written in the dual basis of $\theta_1,\ldots,\theta_n$. 
The associated involutions are the conjugations by 
\[ \sigma_1=\exp(\pi\sqrt{-1}v_1)\]
and
\[ \sigma_j=\exp(2\pi\sqrt{-1}v_j)\]
for $j=2,\ldots,n$. The complex spin representation has weights
$\frac12(\pm\theta_1\pm\cdots\pm\theta_n)$. The center of $\Spin{2n+1}$
is $\{\pm1\}$, so each $\sigma_j$ falls into one of two cases:

(i) $\sigma_j^2=-1$. $\Delta_{2n+1}$ is faithful, so the eigenvalues of
$\Delta_{2n+1}(\sigma_j)$ are $\pm\sqrt{-1}$. We have
\[ \Delta_{2n+1}(\sigma_1)=e^{\pi\sqrt{-1}(d\Delta_{2n+1})(v_1)} \]
and 
\begin{equation}\label{1}
\Delta_{2n+1}(\sigma_j)=e^{2\pi\sqrt{-1}(d\Delta_{2n+1})(v_j)}
\end{equation}
for $j=2,\ldots,n$, so this cases happens precisely for
$j$ odd. In this case $V^{\sigma_j}=V^{-\sigma_j}=0$. 

(ii) $\sigma_j^2=1$. The eigenvalues of $\Delta_{2n+1}(\sigma_j)$ are $\pm1$;
they have equal multiplicity, as is not difficult to see.  
By~(\ref{1}), this case happens for $j$ even. Now 
$\dim V^{\sigma_j}=\dim V^{-\sigma_j}=\frac12\dim V=\frac12d$, 
where $d=\dim_{\mathbb R}\Delta_{2n+1}$ is given by~(\ref{d}). 
In order to have a nice involution, we need
$\frac 12d\leq 4+ \dim G/G^{\sigma}\leq 4+n^2+n$, 
since the symmetric subgroup of $\Spin{2n+1}$ of minimal dimension 
is $\Spin{n+1}\times\Spin n$. We just check that this inequality 
does not hold for $n\geq6$. In particular, 
 $\Delta_{13}$, $\Delta_{15}$
and~$\Delta_{17}$ have no boundary. 

We next consider arbitrary representations of $\Spin7$,
non-standard examples. 
The only ones within the dimension bound are $\Lambda^3\R^7$
of dimension $35$ and
\begin{equation}\label{101}
 \Biii1{}1 
\end{equation}
of dimension $48$. In order to exclude there
representations, we list the involutions of $\Spin7$. The
vertices of the Cartan polyhedron are
\begin{eqnarray*}
  v_1&=&(1,0,0),\\
  v_2&=&(\x,\x,0),\\
  v_3&=&(\x,\x,\x),
\end{eqnarray*}
and the associated involutions are the conjugations by
\[ \sigma_1=\exp(\pi\sqrt{-1}v_1)\]
and
\[ \sigma_j=\exp(2\pi\sqrt{-1}v_j)\]
for $j=2$, $3$, yielding the symmetric spaces
\[ \frac{\SO7}{\SO2\times\SO5},\  \frac{\SO7}{\SO4\times\SO3},\  \frac{\SO7}{\SO6}, \]
respectively. 

Consider first $\Lambda^3\R^7$. In this case the 
center of $\Spin7$ acts trivially. 
We have $\dim V^{\sigma_1}=\binom53+5=15$
and $\codim V^{\sigma_1}=20>4+10=14$, so $\sigma_1$ is not
a nice involution. Next,  
$\dim V^{\sigma_3}=\binom62=15$
and $\codim V^{\sigma_3}=20>4+6=10$,
so $\sigma_3$ is not
a nice involution. On the other hand,  
$\dim V^{\sigma_2}=1+3\binom42=19$
and $\codim V^{\sigma_2}=16=4+12=16$. This shows that 
$\sigma_2$ could be a nice involution, but if so, 
it would necessarily 
fix a codimension one stratum of the orbit space of type $S^3$
(cf.~\cite[Remark~3.2]{GKW});
however $\Lambda^3\R^7$ admits no such strata  by~\cite[Corollary~13.4]{S}. 

Consider now~(\ref{101}). 
We can determine the fixed point set of
$\sigma_j$ on it by looking at the decomposition
\[ \Biii1{}{}\otimes\Biii{}{}1=\Biii1{}1\oplus\Biii{}{}1. \]
The left hand-side is $\rho_7\otimes\Delta_7$, the tensor product of the 
vector and spin representations. 
Since $\rho_7(\sigma_j)^2=1$
and $\Delta_7(\sigma_j)^2=-1$ for $j=1$, $3$, 
we have  $(\rho_7\otimes\Delta_7(\sigma_j))^2=-1$ and hence the 
fixed point set of $\sigma_j$ on~(\ref{101}) is zero for $j=1$, $3$. 
On the other hand $\rho_7(\sigma_2)^2=\Delta_7(\sigma_2)^2=1$.
The fixed point set of $\sigma_2$ 
on $\rho_7$ has dimension~$3$, and on $\Delta_7$ it has dimension $4$. 
It follows that the fixed point set of $\sigma_2$ on $\rho_7\otimes\Delta_7$ 
is $28$, and that on~(\ref{101}) has 
dimension~$24$. But then its codimension is $24$, which is bigger than
$4+\dim\SO7/(\SO4\times\SO3)=16$, so that $\sigma_2$ is not a nice involution. 
The center of $\Spin7$ is $\{\pm1\}$, and replacing $\sigma_2$ by $-\sigma_2$ 
we have $\rho_7(-\sigma_2)=\rho_7(\sigma_2)$ and 
$\Delta_7(-\sigma_2)=-\Delta_7(\sigma_2)$, so the codimension 
on $-\sigma_2$ on~(\ref{101}) is again $24$ and 
neither $-\sigma_2$ is a nice involution.  
 
We now move to $\Spin9$. The only non-standard representation falling 
into the dimension range is $\Lambda^3\R^9$ of dimension~$84$.  
Here the center of $\Spin9$ acts trivially. 
 The
vertices of the Cartan polyhedron are
\begin{eqnarray*}
  v_1&=&(1,0,0,0),\\
  v_2&=&(\x,\x,0,0),\\
  v_3&=&(\x,\x,\x,0),\\
  v_3&=&(\x,\x,\x,\x),
\end{eqnarray*}
and the associated involutions are rthe conjugations by
\[ \sigma_1=\exp(\pi\sqrt{-1}v_1)\]
and
\[ \sigma_j=\exp(2\pi\sqrt{-1}v_j)\]
for $j=2,\ldots,4$ yielding the symmetric spaces
\[ \frac{\SO9}{\SO2\times\SO7},\  \frac{\SO9}{\SO3\times\SO6},\  \frac{\SO9}{\SO4\times\SO5},\  \frac{\SO9}{\SO8}, \]
respectively. We have $\dim V^{\sigma_j}=42$, $46$, $40$, $28$
and then $\codim V^{\sigma_j}>18$, $22$, $24$, $12$ for $j=1,\ldots,4$, 
respectively, proving that none of these are nice involutions. 

We can now deal with $\Spin{2n+1}$ for $n\geq5$. It suffices
to see that the following representations fall out of the dimension bound
(cf.~\cite[Table~5]{OV}): 
\[ \Bnii{}{}1{}{}{} \]
has dimension $\frac13n(2n-1)(2n+1)$. 
\[ \Bn{}2{}{}{} \]
has dimension $\frac13(n-1)(n+1)(2n+1)(2n+3)$.
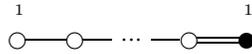
\[ \Bn1{}{}1 \]
has dimension $n\cdot 2^{n+1}$.  
\[ \Bn{}{}{}2 \]
has dimension $\binom{2n+1}n$. 

\subsection{The case of $\SP n$}

Here $\mathcal L_G=48$ for $n=3$, $72$ for $n=4$,
$102$ for $n=5$ and $4n^2$ for $n\geq6$.
The basic examples with non-empty boundary
are the vector representation
$\C^{2n}$, the adjoint representation $\mathfrak{sp}(n)=[S^2\C^{2n}]_{\mathbb R}$, 
and the reduced exterior square $[\Lambda^2_0(\C^{2n})]_{\mathbb R}$.
The special examples are the reduced exterior
powers~$(\Spin4,[\Lambda^4_0\C^8]_{\mathbb R})$ and~$(\Spin6,\Lambda^3_0\C^6)$.

We analyse:
\[ \Cnii{}{}1{}{}{} \]
has dimension $\frac23n(2n+1)(2n-4)$ and is out of the dimension
bound for $n\geq4$; 
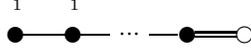
\[ \Cn11{}{} \]
has dimension $\frac{16}3n(n^2-1)$ and is out of the dimension bound for 
$n\geq3$; 
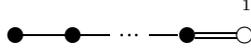
\[ \Cn{}{}{}1 \]
has complex dimension $\frac{2}{n!}(2n+1)2n\cdots(n+3)$,
and it is of quaternionic type if $n$ is odd, and it is of real  
type if $n$ is even, so it is out of the dimension bound for $n\geq5$. 
We get no further examples with non-empty boundary. 

\subsection{The case of $\Spin{2n}$}

Here $\mathcal L_G=4n^2+6n$ for $n\geq4$.
The basic examples are the vector representation
$\R^{2n}$, the adjoint representation $\mathfrak{so}(2n)=\Lambda^2\R^{2n}$, 
and the reduced symmetric square $S^2_0(\R^{2n})$, and 
the special examples are the spin representations
of $\Spin8$, $\Spin{10}$, $\Spin{12}$ and $\Spin{16}$.

We first study the half-spin representations $\Delta_{2n}^\pm$ of $\Spin{2n}$.
\[ \mbox{$\Delta_{2n}$ is of}\quad\begin{array}{ll}\mbox{real type}&\mbox{if $n=4k$,}\\
\mbox{quaternionic type}&\mbox {if $n=4k+2$,}\\
\mbox{complex type}&\mbox {if $n=4k+1$ or $4k+3$.}
\end{array}. \]
Also, $\dim_{\mathbb C}\Delta_{2n}^{\pm\mathbb C}=2^{n-1}$ and
\[ \dim_{\mathbb R}\Delta_{2n}^\pm=\left\{\begin{array}{ll}2^{n-1}&\mbox{if $n=4k$,}\\
  2^n&\mbox {if $n=4k+1$, $4k+2$ or~$4k+3$.}\end{array}.\right. \] 
The only representation 
within the dimension bound which is not a standard
examples is $\Delta_{14}^\pm$.
In the following we apply Borel-de Siebenthal
theory~\cite[Theorem~8.10.8]{Wo}
to describe the involutions of the group and prove that they do not
act as nice involutions on the 
representation. The analysis for $\Spin{14}$ is easily generalized 
to~$\Spin{4k+2}$. 

The roots of $\Spin{14}$ are $\pm\theta_i\pm\theta_j$ for $1\leq i<j\leq7$. 
The system of simple roots
is $\Pi=\{\alpha_1,\ldots,\alpha_7\}$, where $\theta_i=\theta_i-\theta_{i+1}$
for $i=1,\ldots,6$ and $\alpha_7=\theta_6+\theta_7$. The highest root is
$\beta=\theta_1+\theta_2=\alpha_1+2(\alpha_2+\cdots+\alpha_5)+\alpha_6+\alpha_7$. 
Consider the vertices of the Cartan polyhedron
$v_1,\ldots,v_7\in\sqrt{-1}\Lt$, where $\Lt$ is the Lie algebra
of the maximal torus, given by
\[ \alpha_i(v_j)=\left\{\begin{array}{ll}1/m_i,&\mbox{if $i=j$}\\
0,&\mbox{if $i\neq j$}.\end{array}\right. \]
Here the $m_i$ are the coefficients of the highest root.
Then
\begin{eqnarray*}
  v_1&=&(1,0,0,0,0,0,0),\\
  v_2&=&(\x,\x,0,0,0,0,0),\\
  v_3&=&(\x,\x,\x,0,0,0,0),\\
  v_4&=&(\x,\x,\x,\x,0,0,0),\\
  v_5&=&(\x,\x,\x,\x,\x,0,0),\\
  v_6&=&(\x,\x,\x,\x,\x,\x,-\x),\\
  v_7&=&(\x,\x,\x,\x,\x,\x,\x),\\
  \end{eqnarray*}
written in the dual basis of $\theta_1,\ldots,\theta_6$. 
The associated involutions are the conjugations by
\[ \sigma_1=\exp(\pi\sqrt{-1}v_1)\]
for $j=1$, $6$, $7$ and
\[ \sigma_j=\exp(2\pi\sqrt{-1}v_j)\]
for $j=2,\ldots,5$. The half-spin representation $\C^{64}_+$ has real dimension
$128$ and weights $\frac12(\pm\theta_1\pm\cdots\pm\theta_7)$ (mult~$2$),
where the minus sign appears an even number of times.
$\Delta_{14}^\pm$ is faithful and the center of $\Spin{14}$ is 
a cyclic group $\Z_4$; write $\omega$ for a generator. 

For $j=6$, $7$ we have $\sigma_j^2=\pm\omega$ and $\Delta_{14}^+(\sigma_j)$ 
primitive $8$th roots of unity as eigenvalues, so 
$V^{\sigma_j}=V^{z\sigma_j}=0$ for any $z$ in the center of $\Spin{14}$. 

For $j=1$, $3$, $5$, we have $\sigma_j^2=-1$ and $\Delta(\sigma_j)$ 
has $\pm\sqrt{-1}$ as eigenvalues. These are of the same multiplicity.
We can find $z$ in the center of $\Spin{14}$ such that $\Delta_{14}^+(z\sigma_j)$ 
has eigenvalues $\pm1$ of the same multiplicity. 
Now $\mathrm{codim} V^{z\sigma_j}=\frac12\dim V=64$. 
On the other hand the respectively symmetric spaces
$\SO{14}/(\SO2\times\SO{12})$, $\SO{14}/(\SO6\times\SO8)$,
$\SO{14}/(\SO4\times\SO{10})$ have dimensions 
$24$, $48$ and $40$, so $\sigma_j$ is not a nice
involution.  

For $j=2$, $4$, we have $\sigma_j^2=1$ and $\Delta_{14}^+(\sigma_j)$ 
has eigenvalues $\pm1$ of equal multiplicity. Also in this case
$\mathrm{codim} V^{z\sigma_j}=\frac12\dim V=64$. 
On the other hand the respectively symmetric spaces
$\SO{14}/(\SO4\times\SO{10})$, $\SO{14}/(\SO6\times\SO8)$,
have dimensions $40$ and $48$, so $\sigma_j$ is not a nice
involution.  

We next consider arbitrary non-polar representations of $\Spin8$.
The only one within the dimension bound is
\begin{equation}\label{d101}
 \Div1{}1{} 
\end{equation}
of dimension $56$. In order to exclude this
representation, we list the involutions of $\Spin8$. The
vertices of the Cartan polyhedron are
\begin{eqnarray*}
  v_1&=&(1,0,0,0),\\
  v_2&=&(\x,\x,0,0),\\
  v_3&=&(\x,\x,\x,-\x),\\
  v_4&=&(\x,\x,\x,\x),
\end{eqnarray*}
and the associated involutions are the conjugations by
\[ \sigma_j=\exp(\pi\sqrt{-1}v_j)\]
for $j=1$, $3$, $4$, and
\[ \sigma_2=\exp(2\pi\sqrt{-1}v_2),\]
yielding the symmetric spaces
\[ \frac{\SO8}{\SO2\times\SO6} \]
in case $j=1$,
\[ \frac{\SO8}{\U4} \]
in cases $j=3$, $4$, and
\[ \frac{\SO8}{\SO4\times\SO4} \]
in case $j=2$.

We recall that the center of $\Spin8$ is $\Z_2\times\Z_2$,
where the three nontrivial elements correspond each one to a 
generator of the kernel of $\rho_8$, $\Delta_8^+$, $\Delta_8^-$.  

We can determine the fixed point set of $\sigma_j$ (or $z\sigma_j$
for $z$ in the center of $\Spin8$) on
(\ref{d101}) by looking at the decomposition
\[ \Div1{}{}{}\otimes\Div{}{}1{}=\Div1{}1{}\oplus\Div{}{}{}1. \]
The left hand-side is $\rho_8\otimes\Delta_8^+$, the tensor product of the 
vector and half-spin representations.
Now $\rho_8(\sigma_1)=\Delta_8^+(\sigma_4)=\Delta_8^-(\sigma_3)=\mathrm{diag}(-1_2,1_6)$
and $\Delta_8^+(\sigma_1)$, $\rho_8(\sigma_3)$,
$\Delta_8^+(\sigma_3)$, $\rho_8(\sigma_4)$ are complex structures
on $\R^8$. It follows that the eigenvalues of $\sigma_1$ and $\sigma_4$ 
on under $\rho_8\otimes\Delta_8^+$ are $\pm\sqrt{-1}$ (each with multiplicity $32$);
in particular, their fixed point set is zero (and the same is true
for $z\sigma_1$, $z\sigma_4$ for any $z$ in the center). As regards $\sigma_3$,
its eigenvalues under $\rho_8\otimes\Delta_8^+$ are $\pm1$, each with 
multiplicity $32$. It follows that the fixed point set of 
$\sigma_3$ (resp.~$z\sigma_3$ for certain $z$ in the center)
in~(\ref{d101}) is of dimension~$26$ (resp.~$30$), therefore the 
codimension is $26$, which is bigger than the bound
$4+\dim\SO8/\U4=16$ allowed for a nice involution. 
It remains to consider $\sigma_2$.  
It has fixed point set of dimension $4$ on 
both $\rho_8$ and $\Delta_8$. 
It follows that the fixed point set of $\sigma_2$ on~(\ref{d101}) has 
dimension~$28$. But then its codimension is $28$, which is bigger than
$4+\dim\SO8/(\SO4\times\SO4)=20$, so that $\sigma_2$ is not a nice involution
($z\sigma_2$ for $z$ in the center gives the same numbers, so neither 
$z\sigma_2$ is a nice involution). 

We next study the involutions of $\Spin{10}$ and apply the result
to prove that there are no nice involutions for the representation
$\Lambda^3\R^{10}$ of dimension~$120$. The center of $\Spin{10}$ 
acts trivially on this representation.
Here the weights are $\pm\theta_i\pm\theta_j\pm\theta_k$ (mult~$1$)
where $1\leq i<j<k\leq5$, and $\pm\theta_i$ (mult~$4$) where $1\leq i\leq5$. 
The vertices of the Cartan polyhedron are
\begin{eqnarray*}
  v_1&=&(1,0,0,0,0),\\
  v_2&=&(\x,\x,0,0,0),\\
  v_3&=&(\x,\x,\x,0,0),\\
  v_4&=&(\x,\x,\x,\x,-\x),\\
  v_5&=&(\x,\x,\x,\x,\x),
\end{eqnarray*}
and the associated involutions are the conjugations by
\[ \sigma_j=\exp(\pi\sqrt{-1}v_j)\]
for $j=1$, $4$, $5$, and
\[ \sigma_j=\exp(2\pi\sqrt{-1}v_j),\]
for $j=2$, $3$. We see that $V^{\sigma_1}=V^{\sigma_4}=V^{\sigma_5}=0$.
Also $\codim V^{\sigma_2}=64>4+\dim\SO{10}/(\SO4\times\SO6)=28$;
$\codim V^{\sigma_3}=56>4+\dim\SO{10}/(\SO6\times\SO4)=28$.
So none of these are nice involutions. 

We can now deal with $\Spin{2n}$ for $n\geq5$. It suffices
to see that the following representations fall out of the dimension bound: 
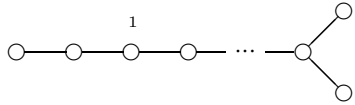
\[ \Dnii{}{}1{}{}{}{} \]
has dimension $\frac23n(n-1)(2n-1)> 4n^2+6n$ for $n\geq6$, and the case $n=5$ was dealt with in the previous paragraph.
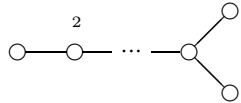
\[ \Dn{}2{}{}{} \]
has dimension $\frac13n(n+1)(2n-3)(2n+1)> 4n^2+6n$ for $n\geq5$.
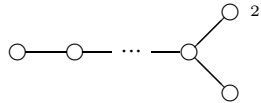
\[ \Dn{}{}{}2{} \]
has dimension $c\binom{2n-1}{n-1}$, where $c=1$ if $n$ is even, and $c=2$
if $n$ is odd; this is bigger than $4n^2+6n$ for $n\geq5$.
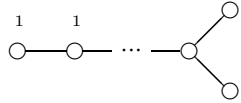
\[ \Dn11{}{}{} \]
has dimension $\frac83n(n-1)(n+1)> 4n^2+6n$ for $n\geq5$.
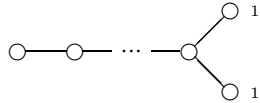
\[ \Dn{}{}{}11 \]
has dimension $\binom{2n}n>4n^2+6n$ for $n\geq5$.
And
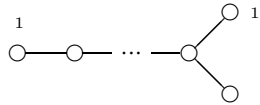
\[ \Dn1{}{}1{} \]
has dimension $c(2n-1)2^{n-1}$, where $c=1$ if $n$ is even, and $c=2$
if $n$ is odd; this is bigger than $4n^2+6n$ for $n\geq5$. We get no further
examples. 

\subsection{Exceptional groups}\label{sec:excep}

Consider $\G$. The standard examples are the polar representations
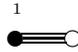
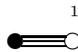
\[ \Gii1{} \]
on $\R^7$
\[ \Gii{}1 \]
and on $\R^{14}$ (adjoint).    
Here  $\mathcal L_G=36$. Since
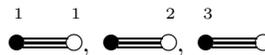
\[ \Gii11,\ \Gii{}2,\ \Gii3{} \]
have dimensions~$64$, $77$ and $77$, respectively,  
there is only one
non-standard representation within the dimension bound, namely,
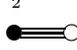
\[ \Gii2{} \]
which is $S^2_0(\R^7)$, of dimension~$27$.

\begin{lem}\label{g2}
  The representation $(\G,\mathrm{S}^2_0\R^7)$ has empty
  boundary in the orbit space.
\end{lem}

\Pf Let $G=\G$ and $V=\mathrm{S}^2_0\R^7$.
It is known (and not difficult to see)
that the principal isotropy group is trivial. Suppose, to the
contrary, that $p$ is a $G$-important point. Then $G_p\approx S^a$
where $a=0$, $1$ or $3$. The first case is ruled out since,
owing to~\cite[Lemma~3.6]{GL},
no $G$-important point may lie on an exceptional orbit.
The last case is excluded
by~\cite[Corollary~13.4]{S}, which implies that the Dynkin index
must be less than~$1$ in case
of a representation of real type a simple Lie group with finite principal isotropy group
and boundary strata of $S^3$-type.

So $G_p\approx S^1$. By the formula in~\cite[Lemma~4.1]{GL},
we can write
\[ \dim V-a-1=\dim G-n+f \]
where $n=\dim N_G(G_p)$ and $f=\dim V^{G_p}$. It gives
$f=n+11\geq13$, as the normalizer of a circle contains a maximal torus.

Let $\alpha_1$, $\alpha_2$ denote the short and long simple roots
of $\G$, respectively. Then the weights of $V$ are
\begin{gather*}
  \pm\alpha_2,\ \pm(3\alpha_1+\alpha_2),\ \pm(3\alpha_1+2\alpha_2),\\
  \pm 2\alpha_1,\ \pm(2\alpha_1+2\alpha_2),\ \pm(4\alpha_1+2\alpha_2)
\end{gather*}
with multiplicity~$1$,
\[ \pm\alpha_1,\ \pm(\alpha_1+\alpha_2),\ \pm(2\alpha_1+\alpha_2) \]
with multiplicity~$2$, and
\[ 0 \]
with multiplicty $3$. The number~$f$ equals the number of
weights that vanish on $\Lg_p$, counted with multiplicity.
It is apparent that $f\leq9$, which is a contradiction. \EPf

\medskip

Consider $\F$.  The standard examples are the polar representations
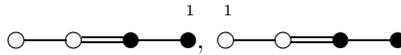
\[ \Fiv{}{}{}1,\ \Fiv1{}{}{} \]
on $\R^{26}$ and on $\R^{52}$ (adjoint), respectively. 
Here  $\mathcal L_G=96$. The next representation
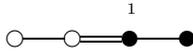
\[ \Fiv{}{}1{} \]
has dimension $273$,
so there are no other representations 
within the dimension bound.

Consider $\E6$.  The standard examples are the polar representations
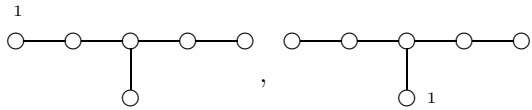
\[ \Evi1{}{}{}{}{},\ \Evi{}1{}{}{}{} \]
on $\C^{27}$ and on $\R^{78}$ (adjoint).
Here  $\mathcal L_G=132$. The next representation
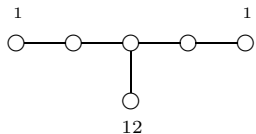
\[ \Evi1{}{}{}{}1 \]
has real dimension $650$,
so there are no other representations 
within the dimension bound.

In the case of $\E7$ we have $\mathcal L_G=222$, and 
the standard examples given by the polar representations
on $\C^{56}$ and $\R^{133}$ (adjoint) are the only ones
within the dimension bound.

In the case of $\E8$ we have $\mathcal L_G=396$, and
the adjoint representation on $\R^{248}$ is the only 
one within the dimension bound.

\section{The reducible case}

According to~\cite[Corollary~7.3]{GKW}, we start by considering
sums of irreducible representations
$(G,V_1)$ and $(G,V_2)$, each of which 
with non-empty boundary in the orbit space.
In view of~\cite[Proposition~7.1]{GKW}, 
$\partial Y_1\neq\varnothing$ or $\partial Y_2\neq\varnothing$,
in the notation there. Hence:

\begin{rmk}\label{finite-pig}
If $H_1$ is finite then none of its elements can act on $V_2$
as a reflection on a hyperplane, because $G$ is connected
so that all of its elements act on $V_2$ by isometries
with determinant $+1$; it follows that $\partial Y_2\neq\varnothing$.
\end{rmk}

\begin{rmk}\label{circle}
If $H_1$ is a torus, then 
there is a circle in $H_1$ which fixes a subspace of real 
codimension two of $V_2$
(cf.~\cite[Lemma~7.1]{GL}).
\end{rmk}

The following two results will rule out many cases. The first one
in fact follows from the second one, but it has a simpler proof,
so we include it. 

\begin{lem}\label{ad0}
Let $G$ be a compact simple Lie group of rank at least two. 
Then $2\mathrm{Ad}$ has empty boundary.
\end{lem}

\Pf In fact $H_1=T^n$ is a maximal torus and acts
on $V_2$ as $n\R \oplus \bigoplus_{\alpha\in\Delta^+}\Lg_{(\alpha)}$ 
where $\Delta^+$ denotes the system of positive roots. 
Suppose, to the contrary, that $2\mathrm{Ad}$ has non-empty boundary.
By Remark~\ref{circle} there is a circle $K$ in $T^n$ which acts on $V_2$ 
with a fixed subspace of real codimension two.
This means $\alpha|_{\mathfrak k}=0$ 
($\Lk$ is the Lie algebra of~$K$) for 
all positive roots $\alpha$ but one, say, $\alpha_0$. Note that~$\Lk$ 
cannot annihilate all simple roots, since these form an integral basis
of~$\Delta^+$. Therefore $\alpha_0$ must be simple. Since the rank 
of $G$ is at least one, there is another simple root $\alpha_1$ 
which is not orthogonal to $\alpha_0$. Now 
$\alpha_0+\alpha_1\in\Delta^+\setminus\{\alpha_0\}$ 
and it does not vanish on $\Lk$, a contradiction. \EPf

\begin{lem}\label{ad}
Let $\rho=(G,V)$ be a real irreducible representation of a  compact 
connect simple Lie group and fix a maximal torus $T^n$ of $G$. 
Then $X=V/T^n$ has empty boundary unless $\rho$ is orbit-equivalent 
to the standard representation of $\SO m$ on $\R^m$ ($m=3$ 
or $m\geq5$). 
\end{lem} 

\Pf Suppose $\partial X\neq\varnothing$. By Remark~\ref{circle},
there is a circle $K$ in $T^n$ whose fixed point subspace~$V^k$ has
real codimension $2$. Let $x$ be a basis of the Lie algebra $\Lk$. 
Choose an ordering of the roots such that $x$ lies in the closure 
of the positive Weyl chamber. 

If $V$ is of complex or quaternionic type
(cf.~\cite[\S4]{GT2}), $V^K$ is the sum of weight spaces whose 
associated weight annihilates~$x$
and $V^K$ has complex codimension~$1$. Then there is a weight $\mu$ 
of multiplicity~$1$ of $V$ such that $\mu(x)\neq0$ and $\lambda(x)=0$ for
every weight $\lambda$ of $V$ with $\lambda\neq\mu$. In particular, 
this shows $V$ cannot be of quaternionic type, because in such a 
case $-\mu$ would be a weight different from $\mu$ and $-\mu(x)\neq0$. 

If $V$ is of real type then the complexification $V^c$ is irreducible,
$(V^c)^K$ is the sum of weight spaces whose 
associated weight annihilates~$x$, and $(V^c)^K$ has complex 
codimension~$2$. There is a weight $\mu$ of multiplicity~$1$ of $V^c$ 
such that $\mu(x)\neq0$ and $\lambda(x)=0$ for
every weight $\lambda$ of $V^c$ with $\lambda\neq\pm\mu$
(here also $-\mu$ is a weight of $V^c$). 

In the complex type case, up to replacing $V$ by $V^*$, which is 
orbit-equivalent to $V$ and has as weights the negatives of the 
weights of $V$, we may assume that $\mu(x)>0$. In the real 
type case, we just rename the weights. In any case, we may assume
that $\mu(x)>0$; we then claim $\mu$ is the highest weight.  
In fact, if $\alpha$ is a positive root then $\alpha(x)\geq0$ 
implies $\mu+\alpha(x)>0$ and $\mu+\alpha\neq\mu$, so 
$\mu+\alpha$ cannot be a weight. 

Continuing, for any positive root $\alpha$, if 
$\langle \mu,\alpha\rangle\neq0$ (here $\langle\cdot,\cdot\rangle$
denotes the Killing form) then $\mu-\alpha$ is a weight, so 
it vanishes on $x$ and hence $\mu(x)=\alpha(x)$.  In particular,
denote by $\alpha_1,\ldots,\alpha_n$ the simple roots, and 
by $\lambda_1,\ldots,\lambda_n$ the fundamental roots. 
Then $\mu=\sum_{i=1}^nm_i\lambda_i$, where 
$m_i=\frac{2\langle\mu,\alpha_i\rangle}{||\alpha_i||^2}$ are non-negative
integers. The above argument shows that if $m_i\neq0$ for some
$i$ then $\mu(x)=\alpha_i(x)$. Let $i_0$ be an index such that 
$m_{i_0}\neq0$. 
Since the Dynkin diagram of $G$ is connected,
for any simple root $\alpha_j$ with $j\neq i_0$, 
we can find a (possible empty) set of 
simple roots $\alpha_{i_1},\ldots,\alpha_{i_p}$ 
such that $\alpha_{i_0}+\alpha_{i_1}+\ldots+\alpha_{i_p}+\alpha_j$ 
is a root~$\beta$. Then 
\[ \langle \mu,\beta\rangle = \langle \mu,\alpha_{i_0}\rangle 
+ \left(\sum_{k=1}^p\langle \mu,\alpha_{i_k}\rangle \right)
+ \langle \mu,\alpha_j\rangle \]
is positive because the first term is positive and others are 
non-negative. Further
\[ \beta(x) = \alpha_{i_0}(x) +
\left(\sum_{k=1}^p \alpha_{i_k}(x) \right)
 + \alpha_j(x), \]
 and $\mu(x)=\beta(x)=\alpha_{i_0}(x)$ implies that
 $\alpha_j(x)=0$ and  $\langle\mu,\alpha_j\rangle=0$. and 
Since $j\neq i_0$ is arbitrary, we have proved that 
$\mu = m_{i_0}\lambda_{i_0}$ and $x=h_{\lambda_{i_0}}$, 
up to multiplication of 
$x$ be a positive scalar, where $h_{\lambda_{i_0}}$
corresponds to $\lambda_{i_0}$ under the musical isomorphism 
induced by the Killing form.  
Now $\mu(x)=\alpha_{i_0}(x)$ implies that 
\[ m_{i_0}\lambda_{i_0}(h_{\lambda_{i_0}})=\alpha_{i_0}(h_{\lambda_{i_0}})=
\langle \alpha_{i_0},\lambda_{i_0}\rangle, \]
so we arrive at the equation
\begin{equation}\label{la}
 2m_{i_0}||\lambda_{i_0}||^2=||\alpha_{i_0}||^2. 
\end{equation}
A quick run of equation~(\ref{la}) 
thorough the fundamental weights of the 
simple Lie algebras (see e.g.~\cite[App.~A]{GT2}) yields as solutions
precisely the
following representations: 
the adjoint representation of $\SU2$, a real form 
of $(\SU4,\Lambda^2\C^4)$, the 
vector representations of $\SO m$ ($m=3$ or $m\geq5$) and the 
half-spin representations $\Delta_8^\pm$. These
are precisely those representations orbit-equivalent to $(\SO m,\R^m)$,
where $m=3$ or $m\geq5$.

Conversely, let us check that $(\SO m,\R^m\oplus\Lambda^2\R^m)$
for $m\geq3$ has non-empty boundary. This follows from~\cite[Proposition~7.1]{GKW}.
In fact $H_1=\SO{m-1}$ acts on $V_2$ as $\R^{m-1}\oplus\Lambda^2\R^{m-1}$, so we can apply
induction and it suffices to see that $(\SO3,\R^3\oplus\Lambda^2\R^3)$
has non-empty boundary, which is clear, as this representation is orbit-equivalent
to $2\R^3$. \EPf

\subsection{The case of $\SU n$}\label{sec:A}
We first consider basic sums, that is, sums such that all summands
are basic examples of irreducible representations with non-empty
boundary, namely, $\C^n$, $\mathrm{Ad}$, $S^2\C^n$, $\Lambda^2\C^n$
($n\geq5$ in the last case). 

Consider $k\C^n$ for $k\geq2$ and $n\geq2$.
The summand $(k-1)\C^n$ has nontrivial
p.i.g.~if and only $k<n$, in which case it is isomorphic to $\SU{n-k+1}$ and
acts on $\C^n$ as $\C^{n-k+1}\oplus2(k-1)\R$ (standard plus trivial), which
has non-empty boundary. By~\cite[Proposition~7.1]{GKW} $k\C^n$
has non-empty boundary if $k<n$.

Next we use again~\cite[Proposition~7.1]{GKW} to show that $n\C^n$
has empty boundary. Let $V_1=(n-1)\C^n$, $V_2=\C^n$.
Then $H_1=\{1\}$, and the action of $H_2=\SU{n-1}$ on $V_1$ is $(n-1)\C^{n-1}$,
up to trivial components, so we proceed by induction and it suffices to show
that $(\SU2,2\C^2)$ has empty boundary. But this is clear, since
the orbit space is the cone over $\Q P^1$. 

Consider $\C^n\oplus\Lambda^2\C^n$ for $n\geq4$. The first summand has p.i.g.\
$\SU{n-1}$ which acts on the second summand as
$\Lambda^2\C^{n-1}\oplus\C^{n-1}$. By ~\cite[Proposition~7.1]{GKW}
and induction, we need only
observe that $\SU3$ on $\C^3\oplus\Lambda^2\C^3=\C^3\oplus\C^{3*}=2\C^3$
has non-empty boundary to deduce that  $\C^n\oplus\Lambda^2\C^n$
has non-empty boundary. 

On the other hand $2\C^n\oplus\Lambda^2\C^n$ for $n\geq4$ has empty boundary
by~\cite[Proposition~7.1]{GKW}. In fact we put $V_1=\C^n$,
$V_2=\C^n\oplus\Lambda^2\C^n$ and check that
both $\partial Y_1=\partial Y_2=\varnothing$. We have $H_2=\{1\}$.
Moreover $H_1=\SU{n-1}$ and it acts on $V_2$ as
$\C^{n-1}\oplus\Lambda^2\C^{n-1}$ plus trivial components, so we
can apply induction and reduce the problem to showing that
$(\SU4,2\C^4\oplus\Lambda^2\C^4)$ has no boundary. This
representation is $2\C^4\oplus2\R^6$; let
$V_1=\C^4$ and $V_2=\C^4\oplus2\R^6$. Then $H_1=\SU3$
acts on $V_2$ as $3\C^3$, which has empty boundary, and
$H_2=\{1\}$. It follows from~\cite[Proposition~7.1]{GKW}
that $2\C^4\oplus2\R^6$ has empty boundary. 

Next we consider $\Lambda^2\C^n\oplus\Lambda^2\C^n$ for $n\geq5$.
We have $H_1=\SU2^m$, where $m=\lfloor \frac n2\rfloor$ and
$(\rho_2(H_1),V_2)$ is
\begin{equation}\label{ll}
 \bigoplus_{i,j=1\atop i<j}^m\C^2_{(i)}\otimes\C^2_{(j)}\oplus m\C \oplus\left[\bigoplus_{i=1}^m\C^2_{(i)}\right], 
\end{equation}
where $\C^2_{(i)}$ is acted by the $i$th-factor of $H_1$ and the term in brackets
exists only if $n$ is odd. The term in brackets clearly has trivial p.i.g..
Further, $\C^2\otimes\C^2$ is orbit-equivalent to the complexification
of $(\SO4,\R^4)$, that is, $\R^4\oplus\R^4$. From this we easily see that 
the first summand in~(\ref{ll}) also has trivial p.i.g.. It follows 
from~\cite[Propositon~7.1]{GKW} that this representation has empty boundary. 

$S^2\C^n$ has finite p.i.g., so it cannot contribute to the
boundary in a sum with another representation,
due to Remark~\ref{finite-pig}.
It follows that $S^2\C^n\oplus S^2\C^n$ for $n\geq2$ has empty boundary. 

Also $\C^n\oplus S^2\C^n$ for $n\geq3$ has empty boundary. Indeed we only
need to show that $\partial Y_2=\varnothing$.
$H_1=\SU{n-1}$ so $V_2|_{\rho_2(H_1)}=\C^{n-1}\oplus  S^2\C^{n-1}$
and we can proceed by induction. Now $(\SU2,\C^2\oplus2\R^3)$ clearly 
has empty boundary. 

We claim $S^2\C^n\oplus\Lambda^2\C^n$ for $n\geq5$ has empty boudary. 
In fact we show $\partial Y_1=\varnothing$. 
$H_2=\SU2^m$, where $m=\lfloor \frac n2\rfloor$ and
$(\rho_1(H_2),V_1)$ is
\begin{equation}\label{ls}
  \bigoplus_{i,j=1\atop i<j}^m\C^2_{(i)}\otimes\C^2_{(j)}\oplus m\C \oplus\left[\bigoplus_{i=1}^mS^2\C^2_{(i)}\right], 
\end{equation}
where $\C^2_{(i)}$ is acted by the $i$th-factor of $H_1$ and the term in brackets
exists only if $n$ is odd.  The term in brackets clearly has trivial p.i.g..
As in the case of~(\ref{ll}), the first summand in~(\ref{ls}) has trivial 
p.i.g., and this proves the claim. 

Next, we show that $\mathrm{Ad}$ cannot be a summand in a representation 
with non-empty boundary together with a basic representation of $\SU n$.
In view of ~\cite[Proposition~7.1]{GKW},
Remark~\ref{finite-pig} and Lemma~\ref{ad}, we 
may suppose that $V_1$ is $\mathrm{Ad}$ and $V_2$ is $\C^n$ or
$\Lambda^2\C^n$, and we only need to show that the 
action of $H_2$ on $V_1=\mathrm{Ad}$  has empty boundary. 
In case $V_2=\C^n$ 
we have $H_2=\SU{n-1}$, $\mathrm{Ad}_{\mathsf{SU}(n)}|_{\mathsf{SU}(n-1)}= 
\mathrm{Ad}_{\mathsf{SU}(n-1)}\oplus\C^{n-1}\oplus\R$, we can use induction,
and the initial case $(\SU2,\R^3\oplus\C^2)$ has empty boundary directly
from~\cite[Proposition~7.1]{GKW}. In case $V_2=\Lambda^2\C^n$ 
($n\geq5$), 
$H_2=\SU2^m$, where $m=\lfloor \frac n2\rfloor\geq2$, and $(\rho_1(H_2),V_1)$ is 
 \begin{equation}\label{al}
  \bigoplus_{i=1}^m\R^3_{(i)}\oplus\bigoplus_{i,j=1\atop i<j}^m\R^4_{(ij)}\oplus (m-1)\R \oplus\left[\R\oplus\bigoplus_{i=1}^m\C^2_{(i)}\right], 
\end{equation}
where $\R^4_{(ij)}$ is acted by the $i$th- and $j$th-factors, $\R^3_{(i)}$ denotes
the adjoint representation of the $i$th factor of $H_1$ and the term in brackets
exists only if $n$ is odd. This representation has empty boundary
by \cite[Proposition~7.1]{GKW},  
because the first two summands make together a representation
with trivial p.i.g., and the last summand also has trivial p.i.g..

Next we consider special sums, that is, sums such that at least one summand is 
a special representation with non-empty  boundary, that is, one
of $(\SU2,\R^5)$, $(\SU4,\R^6)$, $(\SU4,S^2_0(\R^6))$, $(\SU6,\Lambda^3\C^6)$, 
$(\SU8,[\Lambda^4\C^8]_{\mathbb R})$. Note that the first, third and last
representations in this list
have finite p.i.g., so Remark~\ref{finite-pig} applies to them. 

Suppose $(G,V_1)=(\SU2,\R^5)$. 

$\R^5\oplus\R^5$ and $\R^5\oplus\C^2$ 
have empty boundary by Remark~\ref{finite-pig}.

$\R^5\oplus\R^3$ has empty boundary by Lemma~\ref{ad}. 

Suppose $(G,V_1)=(\SU4,\R^6)$. 

It follows easily
from~\cite[Proposition~7.1]{GKW} that $k\R^6$ has boundary if and only
$k\leq5$ and $k\R^6\oplus\ell\C^4$ has boundary if and only if $k+\ell\leq3$. 

$\R^6\oplus\mathfrak{su}(4)$ has non-empty boundary, because $H_2=T^3$
and $(T^3,\R^6)$ is polar so that $\partial Y_1\neq\varnothing$.
But it has trivial p.i.g., and it is easily seen that 
no further summands yield non-empty boundary.

$\R^6\oplus S^2(\C^4)$ has empty boundary. In fact we need
only see that $\partial Y_2=\varnothing$. $H_1=\SP2$ and
$(\rho_2(H_1), V_2)$ is $\mathrm{Ad}_{\mathsf{Sp}(2)}\oplus\mathrm{Ad}_{\mathsf{Sp}(2)}$
which has no boundary by Lemma~\ref{ad0}.

$\R^6\oplus S^2_0(\R^6)$ has empty boundary. In fact we need
only see that $\partial Y_2=\varnothing$. $H_1=\SP2$ acts as $\SO5$
on $V_2$ yielding $S^2_0(\R^5)\oplus\R^5$. Proceeding inductively, we need to
analyse $(\SO3,S^2_0(\R^3)\oplus\R^3)$, which has no boundary
by Remark~\ref{finite-pig} and Lemma~\ref{ad}. 

Suppose $(G,V_1)=(\SU4,S^2_0(\R^6))$.

$S^2_0(\R^6)\oplus S^2_0(\R^6)$ and $S^2_0(\R^6)\oplus S^2(\C^4)$  
have empty boundary by Remark~\ref{finite-pig}.

$S^2_0(\R^6)\oplus\C^4$ has empty boundary.
In fact $H_2=\SU3$ is a simple Lie group and
acts on $V_1$ as $S^2\C^3\oplus\mathrm{Ad}$,
which has empty boundary by the above.

$S^2_0(\R^6)\oplus\mathfrak{so}(6)$ has empty boundary,
by Lemma~\ref{ad}. 

Suppose $(G,V_1)=(\SU6,\Lambda^3\C^6)$.

The p.i.g. $H_1=T^2$ is the maximal torus in the $\SU3$ diagonally embedded
in $\SU3\times\SU3\subset G$. Writing down the weights of
$\C^6$, $\Lambda^2\C^6$, $S^2\C^6$, $\mathfrak{su}(6)$ and $\Lambda^3\C^6$,
we see
that no circle in $H_1$ fixes a subspace of real codimension~$2$,
hence $\partial Y_2=\varnothing$ if $V_2$ is one of
$\C^6$, $\Lambda^2\C^6$, $S^2\C^6$, $\mathfrak{su}(6)$ or~$\Lambda^3\C^6$.  

If $V_2=S^2\C^6$ then $\partial Y_1=\varnothing$ by Remark~\ref{finite-pig}.
If $V_2=\mathfrak{su}(6)$ then $\partial Y_1=\varnothing$ by
Lemma~\ref{ad}.
If $V_2=\C^6$ then $H_2=\SU5$ acts on $V_1$ as $\Lambda^3\C^5\oplus\Lambda^2\C^5$,
which has no boundary by the case of basic sums, so $\partial Y_1=\varnothing$. 
If $V_2=\Lambda^2\C^6$ then $H_2=\SU2^3$ and its action on $V_1$ decomposes
as
\[ \Lambda^3\C^6|_{\mathsf{SU}(2)^3}=2(\C^2_{(1)}\oplus\C^2_{(2)}\oplus\C^2_{(3)})\oplus(\C^2_{(1)}\otimes\C^2_{(2)}\otimes\C^2_{(3)}), \]
where~$\C^2_{(i)}$ is acted by the $i$th factor of $H_2$. This representation
has empty boundary because both summands have trivial p.i.g., so
$\partial Y_1=\varnothing$. 

Suppose $(G,V_1)=(\SU8,[\Lambda^4\C^8]_{\mathbb R})$. Then $H_1=\Z_2^7$ is finite
and we only need to show $\partial Y_1=\varnothing$ below. 

$[\Lambda^4\C^8]_{\mathbb R}\oplus[\Lambda^4\C^8]_{\mathbb R}$ and
$[\Lambda^4\C^8]_{\mathbb R}\oplus S^2\C^8$ have empty boundary
by Remark~\ref{finite-pig}.

$[\Lambda^4\C^8]_{\mathbb R}\oplus\C^8$ has empty boundary because
$H_1=\SU7$ is a simple Lie group and acts on $V_1$ as
$[\Lambda^3\C^7\oplus\Lambda^3\C^{7*}]_{\mathbb R}=\Lambda^3\C^7$,
which has empty boundary by the classification in the irreducible case,
so $\partial Y_1=\varnothing$.

$[\Lambda^4\C^8]_{\mathbb R}\oplus\mathfrak{su}(8)$ has empty boundary 
by Lemma~\ref{ad}.

$[\Lambda^4\C^8]_{\mathbb R}\oplus\Lambda^2\C^8$ has empty boundary because
$H_2=\SU2^4$ and its action on $\Lambda^4\C^8$ decomposes as 
\[ 2\bigoplus_{i<j}\C^2_{(i)}\otimes\C^2_{(j)}\oplus(\C^2_{(1)}\otimes\C^2_{(2)}\otimes\C^2_{(3)}\otimes\C^2_{(4)})\oplus6\C, \]
where~$\C^2_{(i)}$ is acted by the $i$th factor of $H_2$. Therefore
$H_2=\SU2^4$ acts on $V_1=[\Lambda^4\C^8]_{\mathbb R}$ as
\[ \left(2\bigoplus_{i<j}\R^4_{(ij)}\right)\oplus (\R^4_{(12)}\otimes_{\mathbb R}\R^4_{(34)})
\oplus 6\R \]
where $\R^4_{(ij)}$ is acted by the $\SO4$ induced by the $i$th and $j$th
$\SU2$-factors of $H_2$. The first two summands have trivial p.i.g. each,
so that $\partial Y_1=\varnothing$. 

\subsection{The case of $\Spin{2n+1}$}\label{sec:B}
We first consider sums of basic irreducible
representations with non-empty boundary,
those are $\R^{2n+1}$, $\Lambda^2\R^{2n+1}=\mathrm{Ad}$ and $S^2_0(\R^{2n+1})$. 

As in subsection~\ref{sec:A}, $k\R^{2n+1}$ has non-empty boundary precisely for $k\leq 2n$.

$\R^{2n+1}\oplus\Lambda^2\R^{2n+1}$ has non-empty
boundary by Lemma~\ref{ad}.

On the other hand, $2\R^{2n+1}\oplus\Lambda^2\R^{2n+1}$ has empty
boundary. In fact we put $V_1=\R^{2n+1}$, $V_2=\R^{2n+1}\oplus\Lambda^2\R^{2n+1}$.
Note that $H_2=\{1\}$ and $H_1=\Spin{2n}$, which acts on
$V_2$ as $2\R^{2n}\oplus\Lambda^2\R^{2n}$ plus trivial components, so we
can apply induction and reduce the problem to showing that $\SO3$ acting 
on~$2\R^3\oplus\Lambda^2\R^3=2\R^3\oplus\R^{3*}=3\R^3$ has empty boundary,
and this has already been noted. 

$S_2^0\R^{2n+1}\oplus V_2$ has empty boundary if $V_2=S_2^0\R^{2n+1}$,
$\mathrm{Ad}$ or $\R^{2n+1}$. In fact the first case  
follows from Remark~\ref{finite-pig}. In the other two cases we need
only check that~$\partial Y_1=\varnothing$. If $V_2=\mathrm{Ad}$ 
the result follows from Lemma~\ref{ad}. 
If $V_2=\R^{2n+1}$ then $H_2=\Spin{2n}$ acts on $V_1$ as
$S_0^2\R^{2n}\oplus\R^{2n}$ and we proceed by induction to 
analyse~$(\SO3,\R^5\oplus\R^3)$. It was shown in subsection~\ref{sec:A}
that this representation has empty boundary.
  
Next we consider special sums, namely, those involving 
one of the spin representations $(\Spin5=\SP2,\Q^2)$,
$(\Spin7,\R^8)$, 
$(\Spin9,\R^{16})$, $(\Spin{11},\Q^{16})$.
 
Consider $(\SP2,\Q^2)$. 

$2\Q^2$ has non-empty boundary because $H_1=\SP1$ acts 
on $V_2$ as $\Q$ plus trivial components, and this has non-empty 
boundary, but $3\Q^2$ has empty boundary, as is easily seen. 

$\Q^2\oplus\R^5$ has non-empty boundary, because $H_2=\SP1\times\SP1$ 
acts on $V_1$ as $\Q\oplus\Q$ which is polar and this has non-empty boundary. 
This representation has trivial p.i.g., and it is easily seen
thatthis representation
cannot be a summand of another representation with non-empty boundary.

$\Q^2\oplus\Lambda^2\R^5$ has empty boundary. In fact $\partial Y_1=\varnothing$
by Lemma~\ref{ad}.
Further $H_1=\SP1$ acts on $V_2$ as $\R^3\oplus\C^2\oplus3\R$, hence with
empty boundary by the classification in the case of $\SU2$. 

$\Q^2\oplus S^2_0\R^5$ has empty boundary. In fact $H_2=\{1\}$, and $H_1$
acts on $V_2$ as $3\R^3\oplus\C^2\oplus\R$, hence with
empty boundary by the classification in case the case of $\SU2$. 

Consider $\Delta_7=(\Spin7,\R^8)$.

It follows easily from~\cite[Proposition~7.1]{GKW} that
$k\R^7\oplus\ell\R^8$ has non-empty boundary precisely for $k+\ell\leq4$. 
For instance, consider $4\R^8$, which we write as
$V_1\oplus V_2$ for $V_1=V_2=2\R^8$. Then $H_1=\SU3$
acts on $V_2$ as $2\C^3$ plus trivial components, which has non-empty boundary. 

$\R^8\oplus S^2_0\R^7$ has empty boundary. In fact $H_2=\{1\}$, and
$H_1=\G$ acts on $V_2$ with empty boundary by subsection~\ref{sec:excep}.

$\R^8\oplus\Lambda^2\R^7$ has empty boundary. In fact $H_1=\G$ acts on
$V_2$ by $\mathrm{Ad}_{\mathsf{G}_2}\oplus\R^7$, which has empty boundary.
Further, $\partial Y_1=\varnothing$ by Lemma~\ref{ad}.

Consider $\Delta_9=(\Spin9,\R^{16})$. 
The p.i.g. of $\Spin9$ on $\R^{16}$ is a $\Spin7$-subgroup which acts 
on $\R^{16}$ as $\R\oplus\R^7\oplus\R^8$. 

$3\R^{16}$ has non-empty boundary. 
Indeed write $V_1=2\R^{16}$ and $V_2=\R^{16}$.
Then $H_1=\SU3$ acts on $V_2$ as $2\C^3\oplus4\R$ 
which as non-empty boundary and trivial p.i.g.
Similarly one sees that 
$\R^{16}\oplus 4\R^9$ and $2\R^{16}\oplus 2\R^9$ 
have non-empty boundary and trivial p.i.g..
On the other hand, neither one of those
representations, nor $3\R^{16}$ can be a proper summand of 
a representation with non-empty boundary, as is easily checked. 

$\R^{16}\oplus S^2_0(\R^9)$ has empty boundary. In fact $H_2$ is finite, and
$H_1=\Spin7$ acts on $V_2$ with empty boundary, because (cf.~\cite[p.~134]{Fr})
\[ S^2(\rho_9)|_{\mathsf{Spin}(7)}=S^2(\Delta_7\oplus\R)=S^2(\Delta_7)\oplus\Delta_7\oplus\R
=\Lambda^3\rho_7\oplus\Delta_7\oplus2\R, \]
so $V_2|_{H_1}=\Lambda^3\rho_7\oplus\Delta_7\oplus\R$

$\R^{16}\oplus\Lambda^2\R^9$ has empty boundary. In fact 
$H_1=\Spin7$ acts on $V_2$ as $\Lambda^2\rho_7\oplus\rho_7\oplus\Delta_7$,
which has empty boundary. Further  
$\partial Y_1=\varnothing$ by Lemma~\ref{ad}.  

Consider $(\Spin{11},\Q^{16})$. It has trivial p.i.g..

$2\Q^{16}$ and $\Q^{16}\oplus S^2_0(\R^{11})$ have trivial p.i.g.\ by 
Remark~\ref{finite-pig}. 

$\Q^{16}\oplus\R^{11}$ has empty boundary. In fact $H_2=\Spin{10}$ acts 
on $V_1$ as $\C^{16}_+\oplus\C^{16}_-$, and this has trivial boundary, as we shall see
in subsection~\ref{sec:D}.

$\Q^{16}\oplus\Lambda^2\R^{11}$ has empty boundary by Lemma~\ref{ad}. 

\subsection{The case of $\SP n$}

We first consider sums of basic irreducible
representations with non-empty boundary,
those are $\C^{2n}$, $[S^2\C^{2n}]_{\mathbb R}=\mathrm{Ad}$ and a real form
of~$\Lambda^2_0(\C^{2n})=\Lambda^2\C^{2n}\ominus\C$.

As in subsection~\ref{sec:A}, $k\C^{2n}$ has non-empty boundary precisely for $k\leq n$.

We claim $\C^{2n}\oplus[\Lambda^2_0(\C^{2n})]_{\mathbb R}$ has non-empty boundary.
In fact $H_2=\SP1^n$ acts on $V_1$ polarly and hence with non-empty boundary.
It also has trivial p.i.g.. Now $2\C^{2n}\oplus[\Lambda^2_0(\C^{2n})]_{\mathbb R}$ has
empty boundary. In fact we put $V_1=\C^{2n}$,
$V_2=\C^{2n}\oplus[\Lambda^2_0(\C^{2n})]_{\mathbb R}$. Then
$H_2=\{1\}$, and $H_1=\SP{n-1}$ and
\[ V_2|_{H_1}= [\Lambda^2_0\C^{2n-2}]_{\mathbb R} \oplus2\C^{2n-2}\oplus\R, \]
so we can apply induction and reduce the problem to showing that
$\SP2=\Spin5$ acting on
$[\Lambda^2_0\C^4]_{\mathbb R}\oplus 2\C^4=\R^5\oplus2\Q^2$ has empty boundary,
which was proved in subsection~\ref{sec:B}. 

$2[\Lambda^2_0(\C^{2n})]_{\mathbb R}$
has non-empty boundary for $n=3$ and empty boundary for $n\geq4$
In fact $H_1=\SP1^n$ acts on $V_2$ as
\[ \bigoplus_{i<j}\R^4_{(ij)}\oplus(n-1)\R, \]
where $\R^4_{(ij)}$ is acted by the $i$th and $j$th factors of $H_1$. 
If $n=3$ then the p.i.g.\ of $\R^4_{(12)}\oplus\R^4_{(13)}$
is an $\SP1$ diagonally embedded in $H_1$ which acts on $\R^4_{(23)}$
as $\R^3\oplus\R$; this is polar and hence has boundary;
note that it has trivial p.i.g..
On the other other hand, if $n\geq4$ then the p.i.g.\ of
$\R^4_{(12)}\oplus\cdots\oplus\R^4_{(1n)}$ is again 
an $\SP1$ diagonally embedded in $H_1$ which acts on 
the remaining summands as $k\R^3$ plus trivial components, where $k\geq3$,
and this has no boundary. 

Similarly one checks that 
$2[\Lambda_0^2(\C^6)]_{\mathbb R}\oplus\C^6$ 
and $3[\Lambda_0^2(\C^6)]_{\mathbb R}$ have empty boundary. 

$\C^{2n}\oplus[S^2\C^{2n}]_{\mathbb R}$ has empty boundary. 
In fact $\partial Y_1=\varnothing$ by Lemma~\ref{ad}.
Further,
$H_1=\SP{n-1}$ acts on $V_2$ as 
\[ [S^2\C^{2n-2}]_{\mathbb R}\oplus\C^{2n-2}\oplus3\R, \]
so we can use induction and reduce the problem to checking that 
$(\SP1,\R^3\oplus\C^2)$ has no boundary, which was done in subsection~\ref{sec:A}. 

$[S^2\C^{2n}]_{\mathbb R}\oplus[\Lambda^2_0(\C^{2n})]_{\mathbb R}$ for $n\geq3$
has empty boundary. 
In fact $\partial Y_2=\varnothing$ by Lemma~\ref{ad}. 
Moreover $H_2=\SP1^n$ acts on $V_1$ as
\[ \bigoplus_i \R^3_{(i)}\oplus\bigoplus_{i<j}\R^4_{(ij)}. \]
This representation has empty boundary, because the p.i.g.\ of the
first summand is a maximal torus $T^n$ 
of $H_2$ which acts on the second summand, of dimension $4\binom n2$,
with no fixed directions, hence cannot have a boundary; further,
the p.i.g.\ of the second summand is trivial for $n\geq4$ and is a circle
diagonally embedded into $H_2$ if $n=3$, and in any case its action
on the first summand has no boundary.

Next we consider special sums, namely, those involving one of
$(\SP3,\Lambda^3_0\C^6)$, $(\SP4,[\Lambda^4_0\C^8]_{\mathbb R})$. 
These representations have finite p.i.g., therefore no multiple of them
can have non-empty boundary. 

$\C^3\oplus\Lambda^3_0(\C^6)$ has empty boundary. 
Indeed $H_1=\SP2$ acts on $V_2$ as
$4\R^5\oplus\C^4$, which has no boundary. 

$\mathrm{Ad}\oplus\Lambda^3_0(\C^6)$ has empty boundary
by Lemma~\ref{ad}. 

$[\Lambda_0^2(\C^6)]_{\mathbb R}\oplus\Lambda^3_0(\C^6)$ has empty boundary. 
Indeed $H_1=\SP1^3$ acts on $V_2$ as
\[ \left(\C^2_{(1)}\oplus\C^2_{(2)}\oplus\C^2_{(3)}\right)\oplus\left(\C^2_{(1)}\otimes\C^2_{(2)}\otimes\C^2_{(3)}\right). \]
This representation has empty boundary since the action on each term
in parenthesis has trivial p.i.g..

$\C^8\oplus[\Lambda^4_0(\C^8)]_{\mathbb R}$ has empty boundary. 
Indeed $H_1=\SP3$ acts on $V_2$ as
$[\Lambda^2_0\C^6]_{\mathbb R}\oplus\Lambda^3_0\C^6$, which has no boundary by the 
above.

$\mathrm{Ad}\oplus[\Lambda^4_0(\C^8)]_{\mathbb R}$ has empty boundary
by Lemma~\ref{ad}. 

$[\Lambda_0^2(\C^6)]_{\mathbb R}\oplus[\Lambda^4_0(\C^8)]_{\mathbb R}$ 
has empty boundary. 
Indeed $H_1=\SP1^4$ acts on $[\Lambda^4_0(\C^8)]_{\mathbb R}$ as
\[ \left(\C^2_{(1)}\otimes\C^2_{(2)}\otimes\C^2_{(3)}\otimes\C^2_{(4)}\right)\oplus
\bigoplus_{i<j}\C^2_{(i)}\otimes\C^2_{(j)}\oplus2\C \]
and hence on $V_2$ as
\[ \left(\bigoplus_{i<j}\R^4_{(ij)}\right)\oplus (\R^4_{(12)}\otimes_{\mathbb R}\R^4_{(34)})
\oplus 2\R. \]
This representation has empty boundary since the action on each term
in parenthesis has trivial p.i.g..

\subsection{The case of $\Spin{2n}$}\label{sec:D}
We first consider sums of basic irreducible
representations with non-empty boundary,
those are $\R^{2n}$, $\Lambda^2\R^{2n}=\mathrm{Ad}$ and $S^2_0(\R^{2n})$. 

As in subsection~\ref{sec:B}, $k\R^{2n}$ has non-empty boundary precisely for $k\leq 2n-1$;  $\R^{2n}\oplus\Lambda^2\R^{2n}$ has non-empty
boundary; $2\R^{2n}\oplus\Lambda^2\R^{2n}$ 
has empty boundary, as well as
$S_2^0\R^{2n}\oplus V_2$ has empty boundary if $V_2=S_2^0\R^{2n}$,
$\mathrm{Ad}$ or $\R^{2n}$. 

Next we consider special sums, namely, those involving 
one of the spin representations $(\Spin8,\R^8_\pm)$,
$(\Spin{10},\C^{16}_\pm)$, 
$(\Spin{12},\Q^{16}_\pm)$, $(\Spin{16},\R^{128})$,
as well as $(\Spin8,S^2_0\R^8_\pm)$ (note that $\Lambda^2\rho_8=\Lambda^2\Delta_8^\pm)$.

Consider $\Delta_8^\pm=(\Spin8,\R^8_\pm)$. 

Up to an automorphism of $\Spin8$, a representation with two inequivalent
$8$-dimensional 
summands is $\rho_8\oplus\Delta_8^+$. This is polar with p.i.g. given by $\G$.
Since the restriction of 
each of $\rho_8$, $\Delta_8^\pm$ to $\G$ is the same $\R^7\oplus\R$, it does
not matter which $8$-dimensional summands we add to $\rho_8\oplus\Delta_8^+$,
and we obtain that $k\rho_8\oplus\ell\Delta_8^+\oplus m\Delta_8^-$ has non-empty 
boundary for $k+\ell+m\leq5$. If we add a sixth summand, we get empty boundary,
because the p.i.g.\ of one summand is $\Spin7$ and acts on the other 
five summands with empty boundary according to subsection~\ref{sec:B}. 

$\rho_8\oplus S^2_0(\Delta^\pm_8)$ has empty boundary. In fact $H_2=\{1\}$, and
$H_1=\Spin7$ acts on $V_2$ as $S^2_0(\Delta_7)=\Lambda^3(\rho_7)$, which has 
empty boundary. 

Consider $\Delta_{10}^\pm=(\Spin{10},\C^{16}_\pm)$. 

$\Delta_{10}^+\oplus\Delta_{10}^\pm$ 
has empty boundary. In fact $H_1=\SU4$ acts on $V_2$ 
as $2\C^4\oplus2\R^6\oplus4\R$, as we can see from the restricted 
root decomposition of the symmetric space $\E6/(\Spin{10}\U1)$,
and this is not polar. 

We also see easily that $\Delta_{10}^+\oplus3\rho_{10}$ has non-empty 
boundary and  $\Delta_{10}^+\oplus4\rho_{10}$ has empty boundary. 

$\Delta_{10}^+\oplus\Lambda^2(\rho_{10})$ has empty boundary. 
Indeed $H_1=\SU4$ acts on $V_2$ as $\Lambda^2\R^6\oplus2\R^6\oplus2\C^4\oplus2\R$, 
which has empty boundary. Further $\partial Y_1=\varnothing$ 
by Lemma~\ref{ad}. 

$\Delta_{10}^+\oplus S^2(\rho_{10})$ has empty boundary. 
Indeed $H_1=\SU4$ acts on $V_2$ as 
$\Lambda^2\R^6\oplus S^2\C^4\oplus2\C^4\oplus4\R$, 
which has empty boundary. Further $H_2=\{1\}$.

Consider $\Delta_{12}^\pm=(\Spin{12},\Q^{16}_\pm)$. It has p.i.g.\
$\SU2^3\subset\SU6\subset\Spin{12}$. 

$\Delta_{12}^+\oplus\Delta_{12}^+$ 
has empty boundary. In fact $H_1=\SU2^3$ acts on $V_2$ 
as
\[ 16\R\oplus4\left(\R_{(12)}^4\oplus
    \R^4_{(23)}\oplus\R^4_{(31)}\right), \]
and this is easily seen to have empty boundary. 

$\Delta_{12}^+\oplus\Delta_{12}^-$ has empty boundary.
In fact $H_1=\SU2^3$ acts on $V_2$ 
as
\[ 4(\C^2_{(1)}\oplus\C^2_{(2)}\oplus\C^2_{(3)})\oplus
  (\C^2_{(1)}\otimes\C^2_{(2)}\otimes\C^2_{(3)}) \]
and this is easily seen to have empty boundary.

We also see easily that $\Delta_{12}^\pm\oplus\rho_{12}$ has non-empty 
boundary, for $H_1=\SU2^3$ acts on $V_2$ as
$\C^2_{(1)}\oplus\C^2_{(2)}\oplus\C^2_{(3)}$,
which is polar (with trivial p.i.g.) and hence has non-empty boundary.
Further $\Delta_{12}^+\oplus2\rho_{12}$ has empty boundary,
for $H_1=\SU2^3$ acts on $V_2$ as
$2(\C^2_{(1)}\oplus\C^2_{(2)}\oplus\C^2_{(3)})$,
and $H_2=\Spin{10}$ acts on $V_1$ as $\Delta_{10}^+\oplus\Delta_{10}^-$,
and both have empty boundary. 

$\Delta_{12}^+\oplus\Lambda^2(\rho_{12})$ has empty boundary. 
Indeed $H_1=\SU2^3$ acts on $V_2$ as
$9\R\oplus\bigoplus_i\R^3_{(i)}\oplus4\bigoplus_{i<j}\R^4_{(ij)}$,
which has empty boundary. Further $\partial Y_1=\varnothing$ 
by Lemma~\ref{ad}. 

$\Delta_{12}^+\oplus S^2(\rho_{12})$ has empty boundary. 
Indeed $H_1=\SU2^3$ acts on $V_2$ as 
$3\R\oplus3\bigoplus_i\R^3_{(i)}\oplus4\bigoplus_{i<j}\R^4_{(ij)}$, 
which has empty boundary. Further $H_2=\{1\}$.

Consider $\Delta_{16}^\pm=(\Spin{16},\R^{128}_\pm)$. It has
finite p.i.g.. It immediately implies that
$\Delta_{16}^+\oplus\Delta_{16}^\pm$
and $\Delta_{16}^\pm\oplus S^2(\rho_{16})$ have empty boundary. 

We also see easily that $\Delta_{16}^\pm\oplus\rho_{16}$ has empty 
boundary, for $H_2$ acts on $V_1$ as $\Delta_{15}$,
and this has empty boundary.

$\Delta_{16}^+\oplus \Lambda^2(\rho_{16})$ has empty boundary
by Lemma~\ref{ad}.

\subsection{The case of exceptional groups}

Consider $\G$.

We claim $3\R^7$ has non-empty boundary. Indeed the p.i.g. of $2\R^7$
is $\SU2$, which acts on $\R^7$ as $\C^2\oplus3\R$, which is polar
and hence has non-empty boundary. On the other hand,
$4\R^7=2\R^7\oplus2\R^7$ has empty boundary, since
$H_1=\SU2$ acts on $V_2$ as  $2\C^2\oplus6\R$, which has empty boundary.

$\mathfrak{g}_2\oplus\R^7$ has empty boundary. In fact $H_1=T^2$
acts on $V_2$ as the maximal torus of $\SU3$ acts on $\C^3$ (plus a
trivial component) and this has empty boundary (cf.~Lemma~3.4). Further, $H_2=\SU3$
acts on $V_1$ as $\mathfrak{su}(3)\oplus\C^3$, and this has
empty boundary, too. 

Consider $\F$.

$2\R^{26}$  has non-empty boundary.
Indeed the p.i.g. of $\R^{26}$
is $\Spin8$, which acts on $\R^{26}$ as
$2\R\oplus\rho_8\oplus\Delta_8^+\oplus\Delta_8^-$.
This has non-empty aboundary according to the results in
subsection~\ref{sec:D}. On the other hand, $3\R^{26}=2\R^{26}\oplus\R^{26}$
has empty boundary. Indeed $H_1=\SU3$ acts on $V_2$
as $8\R\oplus3\C^3$, which has empty boundary; further
$H_2=\Spin8$ acts on $V_1$ as
$4\R\oplus2\rho_8\oplus2\Delta_8^+\oplus2\Delta_8^-$, which also
has empty boundary. 

$\R^{26}\oplus\mathfrak{f}_4$ has empty boundary.
Indeed $H_1=\Spin8$ acts on $V_2$ as
$\mathfrak{so}(8)\oplus\rho_8\oplus\Delta_8^+\oplus\Delta_8^-$,
which has empty boundary; further, $\partial Y_1=\varnothing$
by Lemma~\ref{ad}. 

Consider $\E6$.

$2\C^{27}$ has empty boundary. In fact $H_1=\Spin8$ acts 
on $\C^{27}$ as $6\R\oplus2(\rho_8\oplus\Delta_8^+\oplus\Delta_8^-)$,
and this has empty boundary. 

$\C^{27}\oplus\mathfrak{e}_6$ has empty boundary. 
In fact $H_1=\Spin8$ acts on $V_2$ as 
$\mathfrak{so}(8)\oplus2(\rho_8\oplus\Delta_8^+\oplus\Delta_8^-)\oplus2\R$,
which has empty boundary; further, $\partial Y_1=\varnothing$
by Lemma~\ref{ad}. 

Consider $\E7$.

$2\Q^{28}$ has empty boundary. In fact $H_1=\Spin8$ acts 
on $\Q^{28}$ as $16\R\oplus4(\rho_8\oplus\Delta_8^+\oplus\Delta_8^-)$,
and this has empty boundary by subsection~\ref{sec:D}. 

$\Q^{28}\oplus\mathfrak{e}_7$ has empty boundary. 
In fact $H_1=\Spin8$ acts on $V_2$ as 
$\mathfrak{so}(8)\oplus4(\rho_8\oplus\Delta_8^+\oplus\Delta_8^-)\oplus9\R$,
which has empty boundary; further, $\partial Y_1=\varnothing$
by Lemma~\ref{ad}. 

In case of $\E8$ there are no cases to check in view of 
Lemma~\ref{ad0}. 

\providecommand{\bysame}{\leavevmode\hbox to3em{\hrulefill}\thinspace}
\providecommand{\MR}{\relax\ifhmode\unskip\space\fi MR }
\providecommand{\MRhref}[2]{%
  \href{http://www.ams.org/mathscinet-getitem?mr=#1}{#2}
}
\providecommand{\href}[2]{#2}


\end{document}